# TWO LIKELIHOOD-BASED SEMIPARAMETRIC ESTIMATION METHODS FOR PANEL COUNT DATA WITH COVARIATES

By Jon A. Wellner[1] and Ying Zhang

*University of Washington and University of Iowa*

We consider estimation in a particular semiparametric regression model for the mean of a counting process with "panel count" data. The basic model assumption is that the conditional mean function of the counting process is of the form $E\{\mathbb{N}(t)|Z\} = \exp(\beta_0^T Z)\Lambda_0(t)$ where $Z$ is a vector of covariates and $\Lambda_0$ is the baseline mean function. The "panel count" observation scheme involves observation of the counting process $\mathbb{N}$ for an individual at a random number $K$ of random time points; both the number and the locations of these time points may differ across individuals.

We study semiparametric maximum pseudo-likelihood and maximum likelihood estimators of the unknown parameters $(\beta_0, \Lambda_0)$ derived on the basis of a nonhomogeneous Poisson process assumption. The pseudo-likelihood estimator is fairly easy to compute, while the maximum likelihood estimator poses more challenges from the computational perspective. We study asymptotic properties of both estimators assuming that the proportional mean model holds, but dropping the Poisson process assumption used to derive the estimators. In particular we establish asymptotic normality for the estimators of the regression parameter $\beta_0$ under appropriate hypotheses. The results show that our estimation procedures are robust in the sense that the estimators converge to the truth regardless of the underlying counting process.

**1. Introduction.** Suppose that $\mathbb{N} = \{\mathbb{N}(t) : t \geq 0\}$ is a univariate counting process. In many applications, it is important to estimate the expected number of events $E\{\mathbb{N}(t)|Z\}$ which will occur by the time $t$, conditionally on a covariate vector $Z$.

Received October 2005; revised December 2006.
[1] Supported in part by NSF Grant DMS-02-03320, NIAID Grant 2R01 AI291968-04 and an NWO Grant to the Vrije Universiteit Amsterdam.

In this paper we consider the proportional mean regression model given by

$$\Lambda(t|Z) \equiv E\{\mathbb{N}(t)|Z\} = e^{\beta_0^T Z}\Lambda_0(t), \tag{1.1}$$

where $\Lambda_0$ is a monotone increasing *baseline mean function*. The parameters of primary interest are $\beta_0$ and $\Lambda_0$.

Suppose that we observe the counting process $\mathbb{N}$ at a random number $K$ of random times $0 \equiv T_{K,0} < T_{K,1} < \cdots < T_{K,K}$. We write $\underline{T}_K \equiv (T_{K,1}, \ldots, T_{K,K})$, and we assume that $(K, \underline{T}_K|Z) \sim G(\cdot|Z)$ is conditionally independent of the counting process $\mathbb{N}$ given the covariate vector $Z$. We further assume that $Z \sim H$ on $\mathbb{R}^d$ with some mild conditions on $H$ for the identifiability of our semiparametric regression model given in Section 3.

The observation for each individual consists of

$$X = (Z, K, \underline{T}_K, \mathbb{N}(T_{K,1}), \ldots, \mathbb{N}(T_{K,K})) \equiv (Z, K, \underline{T}_K, \underline{\mathbb{N}}_K). \tag{1.2}$$

This type of data is referred to as *panel count data*. Throughout this manuscript, we will assume that we observe $n$ i.i.d. copies, $X_1, \ldots, X_n$, of $X$.

Panel count data arise in many fields including demographic studies, industrial reliability and clinical trials; see, for example, Kalbfleisch and Lawless [9], Gaver and O'Muircheataigh [5], Thall and Lachin [16], Thall [15], Sun and Kalbfleisch [13] and Wellner and Zhang [21], where the estimation of either the intensity of event recurrence or the mean function of a counting process with panel count data was studied. Many applications involve covariates whose effects on the underlying counting process are of interest. While there is considerable work on regression modeling for recurrent events based on continuous observations (see, e.g., Lawless and Nadeau [10], Cook, Lawless and Nadeau [3] and Lin, Wei, Yang and Ying [11]), regression analysis with panel count data for counting processes has just started recently. Sun and Wei [14] and Hu, Sun and Wei [6] proposed estimating equation methods, while Zhang [25, 26] proposed a pseudo-likelihood method for studying the proportional mean model (1.1) with panel count data.

To derive useful estimators for this model we will often assume, in addition to (1.1), that the counting process $\mathbb{N}$, conditionally on $Z$, is a non-homogeneous Poisson process. But our general perspective will be to study the estimators and other procedures when the Poisson assumption *may be violated* and we assume *only* that the proportional mean assumption (1.1) holds. Such a program was carried out by Wellner and Zhang [21] for estimation of $\Lambda_0$ without any covariates for this panel count observation model.

The outline of the rest of the paper is as follows. In Section 2 we describe two methods of estimation, namely *maximum pseudo-likelihood* and *maximum likelihood estimators* of $(\beta_0, \Lambda_0)$. The basic picture is that the pseudo-likelihood estimator is computationally relatively straightforward and easy to implement, while the maximum likelihood estimators are considerably more difficult, requiring an iterative algorithm in the computation of the profile likelihood. In Section 3, we state the main asymptotic results: strong consistency, rate of convergence and asymptotic normality of $\hat{\beta}_n^{ps}$ and $\hat{\beta}_n$, for the maximum pseudo-likelihood and maximum likelihood estimators $(\hat{\beta}_n^{ps}, \hat{\Lambda}_n^{ps})$ and $(\hat{\beta}_n, \hat{\Lambda}_n)$ of $(\beta_0, \Lambda_0)$ assuming only the proportional mean structure (1.1), but *not assuming* that $\mathbb{N}$ is a Poisson process. These results are proved in Section 5 by use of tools from empirical process theory. Although pseudo-likelihood methods have been studied in the context of parametric models by Lindsay [12] and Cox and Reid [4], not much seems to be known about their behavior in nonparametric and semiparametric settings such as the one studied here, even assuming that the base model holds. In Section 4 we present the results of simulation studies to demonstrate the robustness of the methods and compare the relative efficiency of the two methods. An application of our methods to a bladder tumor study is presented in this section as well. A general theorem concerning asymptotic normality of semiparametric $M$-estimators and a technical lemma upon which the proofs of our main theorems rely, are stated and proved in Sections 6 and 7, respectively.

## 2. Two likelihood-based semiparametric estimation methods.

*Maximum pseudo-likelihood estimation.* To derive our estimators we assume that conditionally on $Z$, $\mathbb{N}$ is a nonhomogeneous Poisson process with mean function given by (1.1). The pseudo-likelihood method for this model uses the marginal distributions of $\mathbb{N}$, conditional on $Z$,

$$P(\mathbb{N}(t) = k|Z) = \frac{\Lambda(t|Z)^k}{k!} \exp(-\Lambda(t|Z)),$$

and *ignores dependence* between $\mathbb{N}(t_1)$, $\mathbb{N}(t_2)$ to obtain the *log pseudo-likelihood*:

$$l_n^{ps}(\beta, \Lambda) = \sum_{i=1}^n \sum_{j=1}^{K_i} \{\mathbb{N}^{(i)}(T_{K_i,j}^{(i)}) \log \Lambda(T_{K_i,j}^{(i)})$$
$$+ \mathbb{N}^{(i)}(T_{K_i,j}^{(i)})\beta^T Z_i - e^{\beta^T Z_i}\Lambda(T_{K_i,j}^{(i)})\}.$$

Let $\mathcal{R} \subset \mathbb{R}^d$ be a bounded and convex set, and let $\mathcal{F}$ be the class of functions

(2.1) $\quad \mathcal{F} \equiv \{\Lambda : [0, \infty) \to [0, \infty) | \Lambda \text{ is monotone nondecreasing, } \Lambda(0) = 0\}.$

Then the maximum pseudo-likelihood estimator $(\hat{\beta}_n^{ps}, \hat{\Lambda}_n^{ps})$ of $(\beta_0, \Lambda_0)$ is given by $(\hat{\beta}_n^{ps}, \hat{\Lambda}_n^{ps}) \equiv \arg\max_{(\beta, \Lambda) \in \mathcal{R} \times \mathcal{F}} l_n^{ps}(\beta, \Lambda)$. This can be implemented in two steps via the usual profile (pseudo-)likelihood. For each fixed value of $\beta$ we set

$$\hat{\Lambda}_n^{ps}(\cdot, \beta) \equiv \underset{\Lambda \in \mathcal{F}}{\arg\max}\, l_n^{ps}(\beta, \Lambda) \tag{2.2}$$

and define $l_n^{ps,\text{profile}}(\beta) \equiv l_n^{ps}(\beta, \hat{\Lambda}_n^{ps}(\cdot, \beta))$. Then $\hat{\beta}_n^{ps} = \arg\max_{\beta \in \mathcal{R}} l_n^{ps,\text{profile}}(\beta)$ and $\hat{\Lambda}_n^{ps} = \hat{\Lambda}_n^{ps}(\cdot, \hat{\beta}_n^{ps})$. Note that $l_n^{ps}(\beta, \Lambda)$ depends on $\Lambda$ only at the observation time points. By convention, we define our estimator $\hat{\Lambda}_n^{ps}$ to be the one that has jumps only at the observation time points to insure uniqueness.

The optimization problem in (2.2) is easily solved and the details of the solution can be found in Zhang [26].

*Maximum likelihood estimation.* Under the assumption that conditionally on $Z$, $\mathbb{N}$ is a nonhomogeneous Poisson process, the likelihood can be calculated using the (conditional) independence of the increments of $\mathbb{N}$, $\Delta\mathbb{N}(s,t] \equiv \mathbb{N}(t) - \mathbb{N}(s)$, and the Poisson distribution of these increments,

$$P(\Delta\mathbb{N}(s,t] = k | Z) = \frac{[\Delta\Lambda((s,t]|Z)]^k}{k!} \exp(-\Delta\Lambda((s,t]|Z)),$$

to obtain the *log-likelihood*

$$l_n(\beta, \Lambda) = \sum_{i=1}^n \sum_{j=1}^{K_i} \{\Delta\mathbb{N}_{K_ij}^{(i)} \cdot \log \Delta\Lambda_{K_ij} + \Delta\mathbb{N}_{K_ij}^{(i)} \beta^T Z_i - e^{\beta^T Z_i} \Delta\Lambda_{K_ij}\},$$

where

$$\Delta\mathbb{N}_{Kj} \equiv \mathbb{N}(T_{K,j}) - \mathbb{N}(T_{K,j-1}), \qquad j = 1, \ldots, K,$$
$$\Delta\Lambda_{Kj} \equiv \Lambda(T_{K,j}) - \Lambda(T_{K,j-1}), \qquad j = 1, \ldots, K.$$

Then $(\hat{\beta}_n, \hat{\Lambda}_n) \equiv \arg\max_{(\beta,\Lambda) \in \mathcal{R} \times \mathcal{F}} l_n(\beta, \Lambda)$. This maximization can also be carried out in two steps via profile likelihood. For each fixed value of $\beta$ we set $\hat{\Lambda}_n(\cdot, \beta) \equiv \arg\max_{\Lambda \in \mathcal{F}} l_n(\beta, \Lambda)$, and define $l_n^{\text{profile}}(\beta) \equiv l_n(\beta, \hat{\Lambda}_n(\cdot, \beta))$. Then $\hat{\beta}_n = \arg\max_{\beta \in \mathcal{R}} l_n^{\text{profile}}(\beta)$ and $\hat{\Lambda}_n = \hat{\Lambda}_n(\cdot, \hat{\beta}_n)$. Similarly, the estimator $\hat{\Lambda}_n$ is defined to have jumps only at the observation time points. To compute the estimate $(\hat{\beta}_n, \hat{\Lambda}_n)$, we adopt a doubly iterative algorithm to update the estimates alternately. The sketch of the algorithm consists of the following steps:

S1. Choose the initial $\beta^{(0)} = \hat{\beta}_n^{ps}$, the maximum pseudo-likelihood estimator.

S2. For given $\beta^{(p)}$ ($p = 0, 1, 2, \ldots$), the updated estimate of $\Lambda_0$, $\Lambda^{(p)}$, is computed by the modified iterative convex minorant algorithm proposed by Jongbloed [8] on the likelihood $l_n(\beta^{(p)}, \Lambda)$. Initialize this algorithm using $\Lambda^{(p-1)}$ and stop the iteration when

$$\left| \frac{l_n(\beta^{(p)}, \Lambda_{\text{new}}) - l_n(\beta^{(p)}, \Lambda_{\text{current}})}{l_n(\beta^{(p)}, \Lambda_{\text{current}})} \right| \leq \eta.$$

In the very first step, we choose the starting value of $\Lambda$ by interpolating $\hat{\Lambda}_n^{ps}$ linearly between two adjacent jump points to make it monotone increasing and so the likelihood $l_n(\beta, \Lambda)$ is well defined.

S3. For given $\Lambda^{(p)}$, the updated estimate of $\beta$, $\beta^{(p+1)}$, is obtained by optimizing $l_n(\beta, \Lambda^{(p)})$ using the Newton–Raphson method. Initialize the algorithm using $\beta^{(p)}$ and stop the iteration when $\|\beta_{\text{new}} - \beta_{\text{current}}\|_\infty \leq \eta$.

S4. Repeat steps 2 and 3 until the following convergence criterion is satisfied:

$$\left| \frac{l_n(\beta^{(p+1)}, \Lambda^{(p+1)}) - l_n(\beta^{(p)}, \Lambda^{(p)})}{l_n(\beta^{(p)}, \Lambda^{(p)})} \right| \leq \eta.$$

As in the case of pseudo-likelihood studied in Zhang [26], it is easy to verify that for any given monotone nondecreasing function $\Lambda$, the likelihood $l_n(\beta, \Lambda)$ is a concave function of the regression parameter $\beta$ with a negative definite Hessian matrix. Using this fact, we can easily show that the iteration process increases the value of the likelihood, that is, $l_n(\beta^{(p+1)}, \Lambda^{(p+1)}) - l_n(\beta^{(p)}, \Lambda^{(p)}) \geq 0$, for $p = 0, 1, \ldots$.

The iterative algorithms proposed via the profile pseudo-likelihood or the profile likelihood approach converge very well and the convergence does not seem to be affected by the starting point in our simulation experiments described in Section 4. However, this algorithm is not efficient, especially for the maximum likelihood estimation method. It generally needs a considerable number of iterations to achieve the convergence criterion as stated in S4. Meanwhile, computing the profile estimator $\hat{\Lambda}_n$ given in S2 involves the modified iterative convex minorant algorithm which also needs a large number of iterations to converge with the criterion stated in S2. Our simulation experiment with sample size of $n = 100$ shows that computing the maximum likelihood estimator with $\eta = 10^{-10}$ needs about 1800 minutes to converge on a PC (Intel Xeon CPU 2.80 GHz) with the algorithm coded in R. Compared to the profile likelihood algorithm, the profile pseudo-likelihood algorithm is computationally less demanding, since the profile pseudo estimator $\hat{\Lambda}_n^{ps}$ has an explicit solution, as shown in Zhang [26], and hence does not involve any iteration. As result, computing the maximum pseudo-likelihood estimator is much faster than computing the maximum likelihood estimator.

**3. Asymptotic theory: results.** In this section we study the properties of the estimators $(\hat{\beta}_n^{ps}, \hat{\Lambda}_n^{ps})$ and $(\hat{\beta}_n, \hat{\Lambda}_n)$. We establish strong consistency and derive the rate of convergence of both estimators in some $L_2$-metrics related to the observation scheme. We also establish the asymptotic normality of both $\hat{\beta}_n^{ps}$ and $\hat{\beta}_n$ under some mild conditions.

First we give some notation. Let $\mathcal{B}_d$ and $\mathcal{B}$ denote the collection of Borel sets in $\mathbb{R}^d$ and $\mathbb{R}$, respectively, and let $\mathcal{B}_1[0, \tau] = \{B \cap [0, \tau] : B \in \mathcal{B}\}$ and $\mathcal{B}_2[0,\tau] = \mathcal{B}_1[0,\tau] \times \mathcal{B}_1[0,\tau]$. On $([0,\tau], \mathcal{B}_1[0,\tau])$ we define measures $\mu_1$, $\mu_2$, $\nu_1$, $\nu_2$ and $\gamma$ as follows: for $B, B_1, B_2 \in \mathcal{B}_1[0,\tau]$ and $C \in \mathcal{B}_d$, set

$$\nu_1(B \times C) = \int_C \sum_{k=1}^{\infty} P(K=k|Z=z)$$
$$\times \sum_{j=1}^{k} P(T_{k,j} \in B | K=k, Z=z) \, dH(z),$$

$$\nu_2(B_1 \times B_2 \times C) = \int_C \sum_{k=1}^{\infty} P(K=k|Z=z)$$
$$\times \sum_{j=1}^{k} P(T_{k,j-1} \in B_1, T_{k,j} \in B_2 | K=k, Z=z) \, dH(z)$$

and

$$\gamma(B) = \int_{\mathbb{R}^d} \sum_{k=1}^{\infty} P(K=k|Z=z) P(T_{k,k} \in B | K=k, Z=z) \, dH(z).$$

We also define the $L_2$-metrics $d_1(\theta_1, \theta_2)$ and $d_2(\theta_1, \theta_2)$ in the parameter space $\Theta = \mathcal{R} \times \mathcal{F}$ as

$$d_1(\theta_1, \theta_2) = \{|\beta_1 - \beta_2|^2 + \|\Lambda_1 - \Lambda_2\|_{L_2(\mu_1)}^2\}^{1/2},$$
$$d_2(\theta_1, \theta_2) = \{|\beta_1 - \beta_2|^2 + \|\Delta\Lambda_1 - \Delta\Lambda_2\|_{L_2(\mu_2)}^2\}^{1/2},$$

where $\mu_1(B) = \nu_1(B \times \mathbb{R}^d)$ and $\mu_2(B_1 \times B_2) = \nu_2(B_1 \times B_2 \times \mathbb{R}^d)$. To establish consistency, we assume that:

C1. The true parameter $\theta_0 = (\beta_0, \Lambda_0) \in \mathcal{R}^{\circ} \times \mathcal{F}$ where $\mathcal{R}^{\circ}$ is the interior of $\mathcal{R}$.
C2. For all $j = 1, \ldots, K$, $K = 1, 2, \ldots$, the observation times $T_{K,j}$ are random variables taking values in the bounded interval $[0, \tau]$ for some $\tau \in (0, \infty)$. The measure $\mu_l \times H$ on $([0,\tau]^l \times \mathbb{R}^d, \mathcal{B}_l[0,\tau] \times \mathcal{B}_d)$ is absolutely continuous with respect to $\nu_l$ for $l = 1, 2$, and $E(K) < \infty$.
C3. The true baseline mean function $\Lambda_0$ satisfies $\Lambda_0(\tau) \leq M$ for some $M \in (0, \infty)$.
C4. The function $M_0^{ps}$ defined by $M_0^{ps}(X) \equiv \sum_{j=1}^{K} \mathbb{N}_{Kj} \log(\mathbb{N}_{Kj})$ satisfies $PM_0^{ps}(X) < \infty$.
C5. The function $M_0$ defined by $M_0(X) \equiv \sum_{j=1}^{K} \Delta\mathbb{N}_{Kj} \log(\Delta\mathbb{N}_{Kj})$ satisfies $PM_0(X) < \infty$.
C6. $\mathcal{Z} \equiv \text{supp}(H)$, the support of $H$, is a bounded set in $\mathbb{R}^d$. [Thus there exists $z_0 > 0$ such that $P(|Z| \leq z_0) = 1$.]
C7. For all $a \in \mathbb{R}^d$, $a \neq 0$, and $c \in \mathbb{R}$, $P(a^T Z \neq c) > 0$.

Condition C7 is needed together with $\mu_l \times H \ll \nu_l$ from C2 to establish identifiability of the semiparametric model.

THEOREM 3.1. *Suppose that conditions* C1–C7 *hold and the conditional mean structure of the counting process* $\mathbb{N}$ *is given by* (1.1). *Then for every $b < \tau$ for which $\mu_1([b, \tau]) > 0$,*

$$d_1((\hat{\beta}_n^{ps}, \hat{\Lambda}_n^{ps} 1_{[0,b]}), (\beta_0, \Lambda_0 1_{[0,b]})) \to 0 \qquad a.s. \ as \ n \to \infty.$$

*In particular, if $\mu_1(\{\tau\}) > 0$, then*

$$d_1((\hat{\beta}_n^{ps}, \hat{\Lambda}_n^{ps}), (\beta_0, \Lambda_0)) \to 0 \qquad a.s. \ as \ n \to \infty.$$

*Moreover, for every $b < \tau$ for which $\gamma([b, \tau]) > 0$,*

$$d_2((\hat{\beta}_n, \hat{\Lambda}_n 1_{[0,b]}), (\beta_0, \Lambda_0 1_{[0,b]})) \to 0 \qquad a.s. \ as \ n \to \infty.$$

*In particular, if $\gamma(\{\tau\}) > 0$, then*

$$d_2((\hat{\beta}_n, \hat{\Lambda}_n), (\beta_0, \Lambda_0)) \to 0 \qquad a.s. \ as \ n \to \infty.$$

REMARK 3.1. Some condition along the lines of the absolute continuity part of C2 is needed. For example, suppose that $\Lambda_0(t) = t^2$, $\beta_0 = 0$, $\Lambda(t) = t$ and $\beta = 1$. Then if we observe at just one time point $T$ (so $K = 1$ with probability 1), and $T = e^Z$ with probability 1, then $\Lambda_0(T) e^{\beta_0 Z} = \Lambda(T) e^{\beta Z}$ almost surely and the model is not identifiable. C2 holds, in particular, if $(K, T_K)$ is independent of $Z$. The conditions on the measure $\mu_2 \times H$ in C2 and C5 are not needed for proving consistency of $\hat{\theta}_n^{ps} = (\hat{\beta}_n^{ps}, \hat{\Lambda}_n^{ps})$, while the conditions on the measure $\mu_1 \times H$ in C2 and C4 are not needed for proving consistency of $\hat{\theta}_n = (\hat{\beta}_n, \hat{\Lambda}_n)$.

To derive the rate of convergence, we also assume that:

C8. For some interval $O[T] = [\sigma, \tau]$ with $\sigma > 0$ and $\Lambda_0(\sigma) > 0$, $P(\bigcap_{j=1}^{K} \{T_{K,j} \in [\sigma, \tau]\}) = 1$.
C9. $P(K \leq k_0) = 1$ for some $k_0 < \infty$.
C10. For some $v_0 \in (0, \infty)$ the function $Z \mapsto E(e^{v_0 \mathbb{N}(\tau)} | Z)$ is uniformly bounded for $Z \in \mathcal{Z}$.
C11. The observation time points are $s_0$-separated: that is, there exists a constant $s_0 > 0$ such that $P(T_{K,j} - T_{K,j-1} \geq s_0 \ \text{for all} \ j = 1, \ldots, K) = 1$. Furthermore, $\mu_1$ is absolutely continuous with respect to Lebesgue measure $\lambda$ with a derivative $\dot{\mu}_1$ satisfying $\dot{\mu}_1(t) \geq c_0 > 0$ for some positive constant $c_0$.
C12. The true baseline mean function $\Lambda_0$ is differentiable and the derivative has positive and finite lower and upper bounds in the observation interval, that is, there exists a constant $0 < f_0 < \infty$ such that $1/f_0 \leq \Lambda_0'(t) \leq f_0 < \infty$ for $t \in O[T]$.

C13. For some $\eta \in (0,1)$, $a^T \operatorname{Var}(Z|U)a \geq \eta a^T E(ZZ^T|U)a$ a.s. for all $a \in \mathbb{R}^d$, where $(U, Z)$ has distribution $\nu_1/\nu_1(\mathbb{R}^+ \times \mathcal{Z})$.

C14. For some $\eta \in (0,1)$, $a^T \operatorname{Var}(Z|U,V)a \geq \eta a^T E(ZZ^T|U,V)a$ a.s. for all $a \in \mathbb{R}^d$, where $(U, V, Z)$ has distribution $\nu_2/\nu_2(\mathbb{R}^{+2} \times \mathcal{Z})$.

THEOREM 3.2. *In addition to the conditions required for the consistency, suppose* C8, C9, C10 *and* C13 *hold with the constant $v_0$ in C10 satisfying $v_0 \geq 4k_0(1 + \delta_0^{ps})^2$ with $\delta_0^{ps} = \sqrt{c_0 \Lambda_0^3(\sigma)/(24 \cdot 8 f_0)}$ and $\mu_1(\{\tau\}) > 0$. Then*

$$n^{1/3} d_1((\hat{\beta}_n^{ps}, \hat{\Lambda}_n^{ps}), (\beta_0, \Lambda_0)) = O_p(1).$$

*Moreover, if conditions* C11, C12, *and* C14 *hold along with the conditions listed above but with the constant $v_0$ in C10 satisfying $v_0 \geq 4k_0(1 + \delta_0)^2$ with $\delta_0 = \sqrt{c_0 s_0^3/(48 \cdot 8^2 \cdot f_0^4)}$ and $\gamma(\{\tau\}) > 0$, it follows that*

$$n^{1/3} d_2((\hat{\beta}_n, \hat{\Lambda}_n), (\beta_0, \Lambda_0)) = O_p(1).$$

REMARK 3.2. Conditions C8, C9, C10, C11 and C12 are sufficient for validity of Theorem 3.2, but they are probably not necessary. Conditions C9 and C10 are mainly used in deriving the rate of convergence when the counting process $\mathbb{N}$ is allowed to be general [but satisfying the mean model (1.1)]. C8 says that all the observations should fall in a fixed interval in which the mean function is bounded away from zero and C9 indicates that the number of observations is bounded. These conditions are generally true in clinical applications. Condition C10 holds for all $v_0 > 0$, if the counting process is uniformly bounded (which can be justified in many applications) or forms a Poisson process, conditionally on covariates. The first part of C11 requires that two adjacent observation times should be at least $s_0$ apart, an assumption which is very reasonable in practice. The second part of C11 implies that the "total observation measure" $\mu_1$ has a strictly positive intensity (or density). C12 requires that the true baseline mean function should be absolutely continuous with bounded intensity function. While C12 is a reasonable assumption in practice, it may be stronger than necessary. We assume C12 mainly for technical convenience in our proofs.

REMARK 3.3. The metrics $d_1$ and $d_2$ are closely related. Since $\sum_{j=1}^k (a_j - b_j)^2 \leq k^2 \sum_{j=1}^k \{(a_j - a_{j-1}) - (b_j - b_{j-1})\}^2$ (see Wellner and Zhang [20] for a proof), the two metrics are equivalent under C9 and therefore the consistency and rate of convergence results for the maximum likelihood estimator $(\hat{\beta}_n, \hat{\Lambda}_n)$ hold under the metric $d_1$ as well.

REMARK 3.4. Condition C13 can be justified in many applications. By the Markov inequality, it is easy to see that condition C7 implies that $E(ZZ^T)$ is a positive definite matrix. Let $E_1$ and $\text{Var}_1$ denote expectations and variances under the probability measure $\nu_1/\nu_1(\mathbb{R}^+ \times \mathcal{Z})$. If we assume that $\text{Var}_1(Z|U)$ is a positive definite matrix, and we set $\lambda_1 = \max\{\text{eigenvalue}(E_1(ZZ^T|U)\}$ and $\lambda_d^* = \min\{\text{eigenvalue}(\text{Var}_1(Z|U))\}$, then $0 < \lambda_d^* \leq \lambda_1$. Therefore, for any $a \in \mathbb{R}^d$,

$$a^T \text{Var}_1(Z|U)a \geq a^T \lambda_d^* a = \frac{\lambda_d^*}{\lambda_1} a^T \lambda_1 a \geq \frac{\lambda_d^*}{\lambda_1} a^T E_1(ZZ^T|U)a.$$

Thus, condition C13 holds by taking $\eta \leq \lambda_d^*/\lambda_1$. Note that both $\lambda_1$ and $\lambda_d^*$ depend on $U$ in general and the argument here works assuming that this ratio has a positive lower bound uniformly in $U$. We can justify C14 similarly.

Although the overall convergence rate for both the maximum pseudo and likelihood estimators is of the order $n^{-1/3}$, the rate of convergence for the estimators of the regression parameter, as usual, may still be $n^{-1/2}$. Similar to the results of Huang [7] for the Cox model with current status data, we can establish asymptotic normality of both $\hat{\beta}_n^{ps}$ and $\hat{\beta}_n$.

THEOREM 3.3. *Under the same conditions assumed in Theorem 3.2, the estimators $\hat{\beta}_n^{ps}$ and $\hat{\beta}_n$ are asymptotically normal,*

$$(3.1) \qquad \sqrt{n}(\hat{\beta}_n - \beta_0) \xrightarrow{d} Z \sim N_d(0, A^{-1}B(A^{-1})^T)$$

*and*

$$(3.2) \qquad \sqrt{n}(\hat{\beta}_n^{ps} - \beta_0) \xrightarrow{d} Z^{ps} \sim N_d(0, (A^{ps})^{-1}B^{ps}((A^{ps})^{-1})^T),$$

*where*

$$B = E\left\{\sum_{j,j'=1}^{K} C_{j,j'}(Z)[Z - R(T_{K,j}, T_{K,j'})]^{\otimes 2}\right\},$$

$$A = E\left\{\sum_{j=1}^{K} \Delta \Lambda_{0Kj} e^{\beta_0^T Z}[Z - R(K, T_{K,j-1}, T_{K,j})]^{\otimes 2}\right\},$$

$$B^{ps} = E\left\{\sum_{j,j'=1}^{K} C_{j,j'}^{ps}(Z)[Z - R^{ps}(K, T_{K,j})][Z - R^{ps}(K, T_{K,j'})]^T\right\},$$

$$A^{ps} = E\left\{\sum_{j=1}^{K} \Lambda_{0Kj} e^{\beta_0^T Z}[Z - R^{ps}(K, T_{K,j})]^{\otimes 2}\right\},$$

*in which* $R(K, T_{K,j}, T_{K,j'}) \equiv E(Ze^{\beta_0^T Z}|K, T_{K,j}, T_{K,j'})/E(e^{\beta_0^T Z}|K, T_{K,j}, T_{K,j'})$, $R^{ps}(K, T_{K,j}) \equiv E(Ze^{\beta_0^T Z}|K, T_{K,j})/E(e^{\beta_0^T Z}|K, T_{K,j})$, $C_{j,j'}(Z) = \text{Cov}[\Delta \mathbb{N}_{Kj}, \Delta \mathbb{N}_{Kj'}|Z, K, \underline{T}_K]$, $C_{j,j'}^{ps}(Z) = \text{Cov}[\mathbb{N}(T_{Kj}), \mathbb{N}(T_{Kj'})|Z, K, T_{K,j}, T_{K,j'}]$, $\Lambda_{0Kj} = \Lambda_0(T_{K,j})$ *and* $\Delta \Lambda_{0,K,j} = \Lambda_0(T_{K,j}) - \Lambda_0(T_{K,j-1})$.

If the counting process is, conditionally given $Z$, a nonhomogeneous Poisson process with conditional mean function given as specified, then $C_{j,j'}(Z) = \Delta\Lambda_{0Kj} e^{\beta_0^T Z} 1\{j = j'\}$. It follows that $B = A = I(\beta_0)$, the information matrix computed in Wellner, Zhang and Liu [23], and hence $A^{-1}B(A^{-1})^T = I^{-1}(\beta_0)$. This implies that the estimator $\hat{\beta}_n$ under the conditional Poisson process is asymptotically efficient. However, since $C_{j,j'}^{ps}(Z) = e^{\beta_0^T Z} \Lambda_{0K(j \wedge j')}$, $B^{ps} \neq A^{ps}$. This shows that the semiparametric maximum pseudo-likelihood estimator $\hat{\beta}_n^{ps}$ will not be asymptotically efficient under the Poisson assumption.

There is, however, a natural "Poisson regression" model for which the maximum pseudo-likelihood estimator is asymptotically efficient: if we simply assume that the conditional distribution of $(\mathbb{N}(T_{K,1}), \ldots, \mathbb{N}(T_{K,K}))$ given $(K, T_{K,1}, \ldots, T_{K,K}, Z)$ is that of a vector of independent Poisson random variables with means given by $\Lambda(T_{K,j}|Z) = \exp(\beta_0^T Z)\Lambda_0(T_{K,j})$ for $j = 1, \ldots, K$, then

$$C_{j,j'}^{ps}(Z) = \text{Cov}[\mathbb{N}(T_{K,j}), \mathbb{N}(T_{K,j'})|Z, K, T_{K,j}, T_{K,j'}] = \Lambda(T_{K,j}|Z)1\{j = j'\}.$$

Hence $B^{ps} = A^{ps} = I_{\text{PoissRegr}}(\beta_0)$ and $\hat{\beta}^{ps}$ is asymptotically efficient for this alternative model. In practice, this occurs when $(\mathbb{N}(T_{K,1}), \mathbb{N}(T_{K,2}), \ldots, \mathbb{N}(T_{K,K}))$ consist of cluster Poisson count data in which the counts within a cluster are independent.

## 4. Numerical results.

4.1. *Simulation studies.* We generated data using the same schemes as those given in Zhang [26]. Monte Carlo bias, standard deviation and mean squared error of the maximum pseudo-likelihood and maximum likelihood estimates are then compared.

SCENARIO 1. In this scenario, the data is $\{(Z_i, K_i, \underline{T}_{K_i}^{(i)}, \underline{\mathbb{N}}_{K_i}^{(i)}) : i = 1, 2, \ldots, n\}$ with $Z_i = (Z_{i,1}, Z_{i,2}, Z_{i,3})$ where, conditionally on $(Z_i, K_i, \underline{T}_{K_i})$, the counts $\underline{\mathbb{N}}_{K_i}^{(i)}$ were generated from a Poisson process. For each subject, we generate data by the following scheme: $Z_{i,1} \sim \text{Unif}(0,1), Z_{i,2} \sim N(0,1), Z_{i,3} \sim \text{Bernoulli}(0.5)$; $K_i$ is sampled randomly from the discrete set, $\{1, 2, 3, 4, 5, 6\}$. Given $K_i$, $\underline{T}_{K_i}^{(i)} = (T_{K_i,1}^{(i)}, T_{K_i,2}^{(i)}, \ldots, T_{K_i,K_i}^{(i)})$ are the order statistics of $K_i$ random observations generated from $\text{Unif}(1, 10)$ and rounded to the second decimal point to make the observation times possibly tied. The panel counts

$$\underline{\mathbb{N}}_{K_i}^{(i)} = (\mathbb{N}^{(i)}(T_{K_i,1}^{(i)}), \mathbb{N}^{(i)}(T_{K_i,2}^{(i)}), \ldots, \mathbb{N}^{(i)}(T_{K_i,K_i}^{(i)}))$$

are generated from the Poisson process with the conditional mean function given by $\Lambda(t|Z_i) = 2t\exp(\beta_0^T Z_i)$, that is,

$$\mathbb{N}^{(i)}(T_{K_i,j}^{(i)}) - \mathbb{N}^{(i)}(T_{K_i,j-1}^{(i)}) \sim \text{Poisson}\{2(T_{K_i,j}^{(i)} - T_{K_i,j-1}^{(i)})\exp(\beta_0^T Z_i)\},$$

where $\beta_0 = (\beta_1, \beta_2, \beta_3)^T = (-1.0, 0.5, 1.5)^T$.

For this scenario, we can directly calculate the asymptotic covariance matrices given in Theorem 3.3, $\Sigma^{ps} \equiv (A^{ps})^{-1} B^{ps} ((A^{ps})^{-1})^T = (1582/17787) W^{-1}$ and $\Sigma = A^{-1} B (A^{-1})^T = A^{-1} = (1260/19179) W^{-1}$, respectively, where $W = E\{e^{\beta_0^T Z} [Z - E(Z e^{\beta_0^T Z})/E(e^{\beta_0^T Z})]^{\otimes 2}\}$. Since it is difficult to evaluate the matrix $W$ analytically, we calculated it numerically using Mathematica (Wolfram [24]) to obtain the following approximate results for the asymptotic covariance matrices:

$$(4.1) \qquad \Sigma^{ps} \approx \begin{pmatrix} 0.571104 & 0.000000 & 0.000000 \\ 0.000000 & 0.045304 & 0.000000 \\ 0.000000 & 0.000000 & 0.303752 \end{pmatrix}$$

and

$$(4.2) \qquad \Sigma \approx \begin{pmatrix} 0.421848 & 0.000000 & 0.000000 \\ 0.000000 & 0.033464 & 0.000000 \\ 0.000000 & 0.000000 & 0.224368 \end{pmatrix}.$$

We conducted simulation studies with sample sizes of $n = 50$ and $n = 100$, respectively. For each case, the Monte Carlo sample bias, standard deviation and mean squared error for the semiparametric estimators of the regression parameters are reported in Table 1. We also include the asymptotic standard errors obtained from (4.1) and (4.2) in Table 1 to compare with the Monte Carlo sample standard deviations. The results show that the sample bias for both estimators is small, the standard deviation and mean squared error are smaller for the maximum likelihood method compared to the pseudo-likelihood method and the latter decrease as $n^{-1/2}$ and $n^{-1}$, respectively, as sample size increases. Moreover, the standard errors of estimates based on asymptotic theory are all close to the corresponding standard deviations based on the Monte Carlo simulations. All of these provide numerical support for our asymptotic results in Theorem 3.3.

Based on the results of 1000 Monte Carlo samples, we plot the pointwise means and 2.5- and 97.5-percentiles of both estimators of the baseline mean function $\Lambda(t) = 2t$ in Figure 1. It clearly shows that both estimators seem to have negligible bias and the maximum likelihood estimator has smaller variability compared to the maximum pseudo-likelihood estimator. When sample size increases, the variability of both estimators decreases accordingly.

SCENARIO 2. In this scenario, the data is $\{(Z_i, K_i, \underline{T}_{K_i}^{(i)}, \underline{\mathbb{N}}_{K_i}^{(i)}) : i = 1, 2, \ldots, n\}$ with $Z_i = (Z_{i,1}, Z_{i,2}, Z_{i,3})$ and, conditionally on $(Z_i, K_i, \underline{T}_{K_i}^{(i)})$, the counts $\underline{\mathbb{N}}_{K_i}^{(i)}$ were generated from a mixed Poisson process. For each subject, $(Z_i, K_i, \underline{T}_{K_i}^{(i)})$ are generated in exactly the same way as in Scenario 1. The panel counts are, however, generated from a homogeneous Poisson process with a random effect on the intensity: given subject $i$ with covariates $Z_i$ and frailty variable $\alpha_i$ (independent of $Z_i$), the counts are generated from the Poisson process with intensity $(\lambda + \alpha_i) \exp(\beta_0^T Z_i)$, where $\lambda = 2.0$ and $\alpha_i \in \{-0.4, 0, 0.4\}$ with probabilities 0.25, 0.5 and 0.25, respectively.

TABLE 1
*Results of the Monte Carlo simulation studies for the regression parameter estimates based on 1000 repeated samples for data generated from the conditional Poisson process*

|  | $n = 50$ | | $n = 100$ | |
|---|---|---|---|---|
|  | **Pseudo-likelihood** | **Likelihood** | **Pseudo-likelihood** | **Likelihood** |
| Estimate of $\beta_1$ | | | | |
| BIAS | 0.0020 | 0.0018 | 0.0017 | 0.0015 |
| SD | 0.1193 | 0.1019 | 0.0806 | 0.0694 |
| ASE | 0.1069 | 0.0919 | 0.0758 | 0.0649 |
| MSE $\times 10^2$ | 1.4236 | 1.0387 | 0.6499 | 0.4819 |
| Estimate of $\beta_2$ | | | | |
| BIAS | $-0.0003$ | $-0.0016$ | 0.0028 | 0.0022 |
| SD | 0.0349 | 0.0294 | 0.0231 | 0.0193 |
| ASE | 0.0301 | 0.0259 | 0.0213 | 0.0183 |
| MSE $\times 10^2$ | 0.1218 | 0.0867 | 0.0541 | 0.0377 |
| Estimate of $\beta_3$ | | | | |
| BIAS | 0.0023 | 0.0011 | 0.0016 | $-0.0009$ |
| SD | 0.0830 | 0.0712 | 0.0579 | 0.0497 |
| ASE | 0.0779 | 0.0670 | 0.0551 | 0.0474 |
| MSE $\times 10^2$ | 0.6894 | 0.5071 | 0.3355 | 0.2471 |

In this scenario, the counting process given only the covariates is not a Poisson process. However, the conditional mean function of the counting process given the covariates still satisfies (1.1) with $\Lambda_0(t) = 2t$ and thus our proposed methods are expected to be valid for this case as well. The asymptotic variances given in Theorem 3.3 for this scenario are

$$\Sigma^{ps} = (A^{ps})^{-1} B^{ps} ((A^{ps})^{-1})^T = \frac{1582}{17787} W^{-1} + \frac{463.12}{17787} W^{-1} \tilde{W} (W^{-1})^T$$

and

$$\Sigma = A^{-1} B (A^{-1})^T = \frac{1260}{19179} W^{-1} + \frac{7917588}{19179^2} W^{-1} \tilde{W} (W^{-1})^T,$$

respectively, where $\tilde{W} = E\{e^{2\beta_0^T Z}[Z - E(Ze^{\beta_0^T Z})/E(e^{\beta_0^T Z})]^{\otimes 2}\}$. Using Mathematica (Wolfram [24]) to calculate the asymptotic covariance matrices nu-

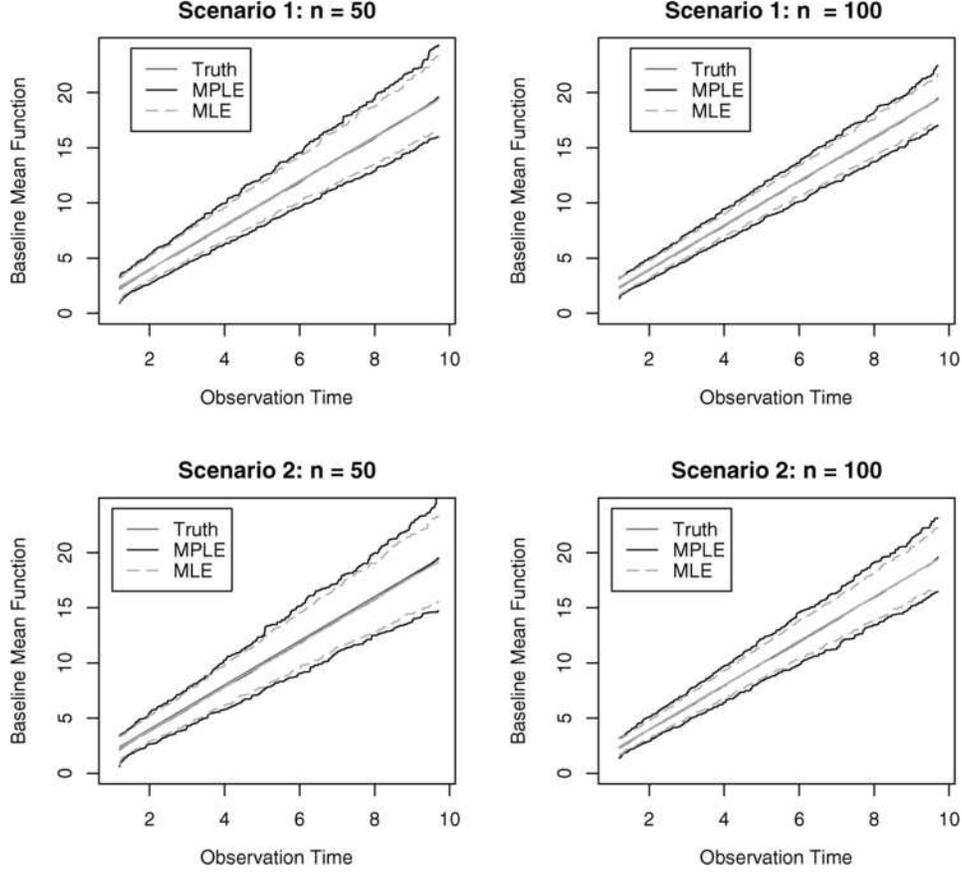

Fig. 1. *The pointwise means, 2.5-percentiles and 97.5-percentiles of both the maximum pseudo-likelihood and likelihood estimators of the baseline mean function under the proportional mean model. MPLE: The maximum pseudo-likelihood estimator; MLE: The maximum likelihood estimator.*

merically yields

$$(4.3) \qquad \Sigma^{ps} \approx \begin{pmatrix} 1.172450 & -0.023852 & -0.043178 \\ -0.023852 & 0.108760 & 0.022975 \\ -0.043178 & 0.022975 & 0.448444 \end{pmatrix}$$

and

$$(4.4) \qquad \Sigma \approx \begin{pmatrix} 0.918986 & -0.019718 & -0.035696 \\ -0.019718 & 0.085924 & 0.018994 \\ -0.035696 & 0.018994 & 0.343985 \end{pmatrix}.$$

As in Scenario 1, we conducted simulation studies with sample sizes of $n = 50$ and $n = 100$, respectively. For each case, the Monte Carlo sample bias, standard deviation and mean squared error for the semiparametric estimators of the regression parameters are computed with 1000 repeated samples. The results are shown in Table 2. In Figure 1, we also plot the pointwise means, 2.5-percentiles and 97.5-percentiles of both estimators of the unconditional baseline mean function $\Lambda_0(t) = 2t$ based on the results obtained from 1000 Monte Carlo samples. We observe the same phenomenon as appeared in Scenario 1: for the regression parameters, both standard deviation and mean squared error using the maximum likelihood method are smaller than those using the pseudo-likelihood method while the bias is relatively small; for the baseline mean function, both estimators have a negligible bias but the maximum likelihood estimator has less variability than the maximum pseudo-likelihood estimator. We also note that the variability results of the semiparametric estimators are relatively larger than their counterpart in Scenario 1. This may be caused by violation of the assumption of a conditional Poisson process given only the covariates. We also include the asymptotic standard errors of the regression parameter estimates based on (4.3) and (4.4) in Table 2. Again the standard errors derived from the asymptotic theory are all close to the standard deviations based on Monte Carlo simulations.

These simulation studies provide numerical support for the statement that the proposed semiparametric estimation methods are robust against the underlying conditional Poisson process assumption. These methods are valid as long as the proportional mean function model (1.1) holds. We have also conducted several analytical analyses to compare the semiparametric efficiency between the maximum pseudo-likelihood and maximum likelihood estimation methods. There is considerable evidence that the maximum likelihood method (based on the Poisson process assumption) is more efficient than the pseudo-likelihood method both on and off the Poisson model, with large differences occurring when $K$ is heavily tailed. The detailed analytical results are presented in Wellner, Zhang and Liu [23].

TABLE 2
*Results of the Monte Carlo simulation studies for the regression parameter estimates based on 1000 repeated samples for data generated from the mixed Poisson process*

|  | $n = 50$ | | $n = 100$ | |
|---|---|---|---|---|
|  | Pseudo-likelihood | Likelihood | Pseudo-likelihood | Likelihood |
| Estimate of $\beta_1$ | | | | |
| BIAS | 0.0038 | 0.0029 | −0.0068 | −0.0072 |
| SD | 0.1556 | 0.1415 | 0.1138 | 0.0993 |
| ASE | 0.1531 | 0.1356 | 0.1083 | 0.0959 |
| MSE×$10^2$ | 2.4226 | 2.0003 | 1.2997 | 0.9912 |
| Estimate of $\beta_2$ | | | | |
| BIAS | −0.0008 | −0.0001 | 0.0012 | 0.0017 |
| SD | 0.0467 | 0.0425 | 0.0318 | 0.0297 |
| ASE | 0.0466 | 0.0415 | 0.0330 | 0.0293 |
| MSE×$10^2$ | 0.2182 | 0.1806 | 0.1013 | 0.0885 |
| Estimate of $\beta_3$ | | | | |
| BIAS | 0.0096 | 0.0099 | 0.0061 | 0.0040 |
| SD | 0.0972 | 0.0888 | 0.0666 | 0.0581 |
| ASE | 0.0947 | 0.0829 | 0.0670 | 0.0587 |
| MSE×$10^2$ | 0.9540 | 0.7983 | 0.4473 | 0.3392 |

4.2. *A real data example.* Using the semiparametric methods proposed in the preceding sections, we analyze the bladder tumor data extracted from Andrews and Herzberg ([1], pages 253–260). This data set comes from a bladder tumor study conducted by the Veterans Administration Cooperative Urological Research (Byar, Blackard and Vacurg [2]). In the study, a randomized clinical trial of three treatments, placebo, pyridoxine pills and thiotepa instillation into the bladder was conducted for patients with superficial bladder tumor when entering the trial. At each follow-up visit, tumors were counted, measured and then removed if observed, and the treatment was continued. The treatment effects, especially the thiotepa instillation, on suppressing the recurrence of bladder tumor have been explored by many authors, for example, Wei, Lin and Weissfeld [19], Sun and Wei [14], Wellner and Zhang [21] and Zhang [26].

In this paper, we study the proportional mean model that has been proposed by Sun and Wei [14] and Zhang [26],

(4.5) $\quad E\{\mathbb{N}(t)|Z\} = \Lambda_0(t)\exp(\beta_1 Z_1 + \beta_2 Z_2 + \beta_3 Z_3 + \beta_4 Z_4),$

where $Z_1$ and $Z_2$ represent the number and size of bladder tumors at the beginning of the trial, and $Z_3$ and $Z_4$ are the indicators for the pyridoxine pill and thiotepa instillation treatments, respectively. We choose $\beta^{(0)} = (0,0,0,0)$ to start our iterative algorithm and $\eta = 10^{-10}$ for the convergence criteria to stop the algorithm. Since the asymptotic variances are difficult to estimate, we adopt the bootstrap procedure to estimate the asymptotic standard error of the semiparametric estimates of the regression parameters. We generated 200 bootstrap samples and calculated the proposed estimators for each sampled data set. The sample standard deviation of the estimates based on these 200 bootstrap samples is used to estimate the asymptotic standard error. The inference based on the bootstrap estimator for asymptotic standard error is given in Table 3. The semiparametric maximum pseudo-likelihood and maximum likelihood estimators of the baseline mean function are plotted in Figure 2.

Both methods yield the same conclusion that the baseline number of tumors (the number of tumors observed when entering the trial) significantly affects the recurrence of the tumor at level 0.05 ($p$-value $= 0.0105$ and 0.0078, resp., for the maximum pseudo-likelihood and maximum likelihood methods), and the thiotepa instillation treatment appears to reduce the recurrence of tumor significantly. ($p$-value $= 0.0186$ and 0.0269, resp., for the maximum pseudo-likelihood and maximum likelihood methods).

In Figure 2, we can see that the maximum likelihood estimator of the baseline mean function is substantially smaller than the maximum pseudo-likelihood estimator, which preserves the phenomenon we have observed in nonparameteric estimation methods for this data set studied in Wellner and Zhang [21].

We also notice that the maximum likelihood method, in contrast to what we have observed through both the simulation and analytical studies, yields larger standard errors compared to the pseudo-likelihood method. Violation of the proportional mean model (4.5) for this data set could be the explanation for this result, since Zhang [27] plotted the nonparameteric pseudo-likelihood estimators of the mean function for each of three treatments and found that the estimators cross over. While plotting the nonparametric estimators of the mean function for the groups defined by covariates is a reasonable first step in an exploration of the validity of the proposed model, it would be preferable to proceed via more quantitative measures, such as appropriate goodness-of-fit statistics. The construction of goodness-of-fit test statistics for regression modeling of panel count data remains an open problem for future research.

All the numerical experiments in this paper were implemented in R. The computing programs are available from the second author.

TABLE 3
*Semiparametric inference for the bladder tumor study based on 200 bootstrap samples from the original data set*

| Variable | Method | $\hat{\beta}$ | $\widehat{se(\hat{\beta})}$ | $\hat{\beta}/\widehat{se(\hat{\beta})}$ | $p$-value |
|---|---|---|---|---|---|
| $Z_1$ | | | | | |
| | Pseudo-likelihood | 0.1446 | 0.0565 | 2.5593 | 0.0105 |
| | Likelihood | 0.2069 | 0.0778 | 2.6594 | 0.0078 |
| $Z_2$ | | | | | |
| | Pseudo-likelihood | $-0.0450$ | 0.0632 | $-0.7120$ | 0.4746 |
| | Likelihood | $-0.0355$ | 0.0861 | $-0.4123$ | 0.6801 |
| $Z_3$ | | | | | |
| | Pseudo-likelihood | 0.1951 | 0.3233 | 0.6035 | 0.5462 |
| | Likelihood | 0.0664 | 0.4310 | 0.1541 | 0.8775 |
| $Z_4$ | | | | | |
| | Pseudo-likelihood | $-0.6881$ | 0.2923 | $-2.3541$ | 0.0186 |
| | Likelihood | $-0.7972$ | 0.3603 | $-2.2126$ | 0.0269 |

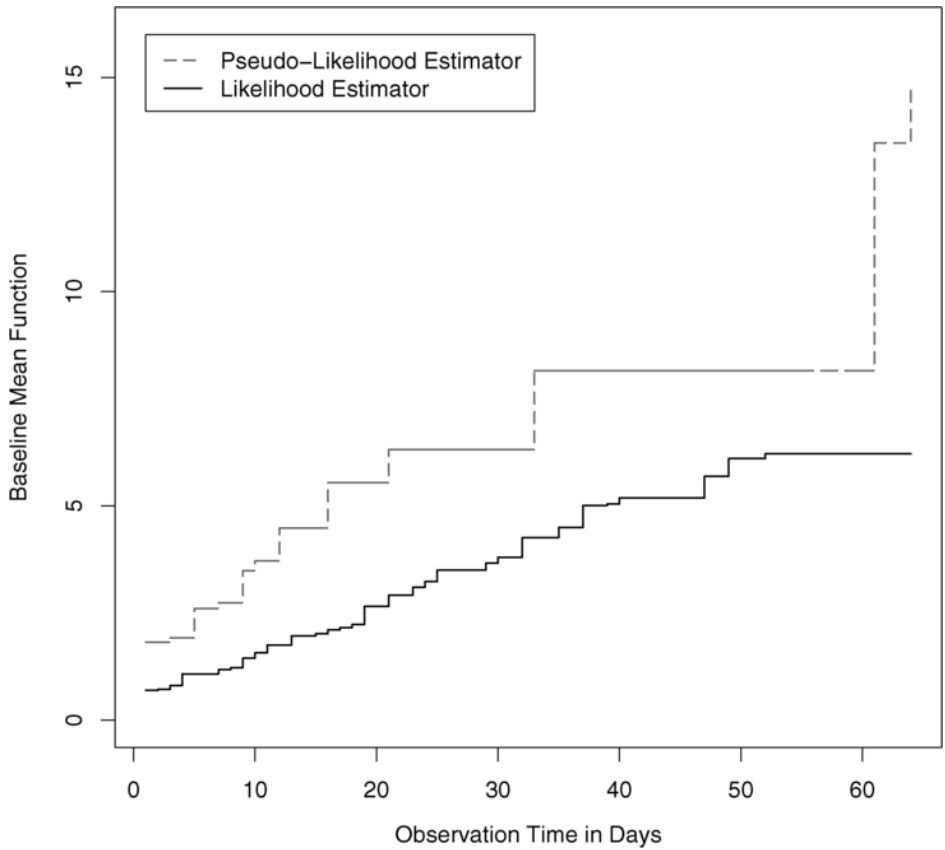

Fig. 2. *The two estimators of the baseline mean function under the proportional mean model for the bladder tumor example.*

**5. Asymptotic theory: proofs.** We use empirical process theory to study the asymptotic properties of the semiparametric maximum pseudo-likelihood and maximum likelihood estimators. The proof of Theorem 3.1 is closely related to the proof of Theorem 4.1 of Wellner and Zhang [21]. The rate of convergence is derived based on the general theorem for the rate of convergence given in Theorem 3.2.5 of van der Vaart and Wellner [18]. The asymptotic normality proofs for both $\hat{\beta}_n^{ps}$ and $\hat{\beta}_n$ are based on the general theorem for $M$-estimation of regression parameters in the presence of a nonparametric nuisance parameter, which is stated (and proved) in Section 6.

PROOF OF THEOREM 3.1. Zhang [26] has given a proof for the first part of the theorem concerning the semiparametric maximum pseudo-likelihood estimator. Unfortunately, his proof of Theorem 1 on pages 47 and 48 is not correct (in particular, the conditions imposed do not suffice for identifiability as claimed). Here we give proofs for both the maximum pseudo-likelihood and maximum likelihood estimators.

We first prove the claims concerning the pseudo-likelihood estimators $(\hat{\beta}_n^{ps}, \hat{\Lambda}_n^{ps})$. Let $\mathbb{M}_n^{ps}(\theta) = n^{-1} l_n(\beta, \Lambda) = \mathbb{P}_n m_\theta^{ps}(X)$ and $\mathbb{M}^{ps}(\theta) = P m_\theta^{ps}(X)$, where

$$m_\theta^{ps}(X) = \sum_{j=1}^K \{\mathbb{N}_{Kj} \log \Lambda_{Kj} + \mathbb{N}_{Kj} \beta^T Z - \Lambda_{Kj} \exp(\beta^T Z)\}.$$

First, we show that $\mathbb{M}^{ps}$ has $\theta_0 = (\beta_0, \Lambda_0)$ as its unique maximizing point. Computing the expectation conditionally on $(Z, K, \underline{T}_K)$ yields

$$\mathbb{M}^{ps}(\theta_0) - \mathbb{M}^{ps}(\theta) = \int \Lambda(u) \exp(\beta^T z) h\left[\frac{\Lambda_0(u) \exp(\beta_0^T z)}{\Lambda(u) \exp(\beta^T z)}\right] d\nu_1(u, z),$$

where $h(x) = x \log(x) - x + 1$. The function $h(x)$ satisfies $h(x) \geq 0$ for $x > 0$ with equality holding only at $x = 1$. Hence $\mathbb{M}^{ps}(\theta_0) \geq \mathbb{M}^{ps}(\theta)$ and $\mathbb{M}^{ps}(\theta_0) = \mathbb{M}^{ps}(\theta)$ if and only if

$$(5.1) \qquad \frac{\Lambda_0(u) \exp(\beta_0^T z)}{\Lambda(u) \exp(\beta^T z)} = 1 \qquad \text{a.e. with respect to } \nu_1.$$

This implies that

$$(5.2) \qquad \beta = \beta_0 \quad \text{and} \quad \Lambda(u) = \Lambda_0(u) \qquad \text{a.e. with respect to } \mu_1$$

by C2 and C7. Here is a proof of this claim: Let

$$f_1(u) = \Lambda(u) - \Lambda_0(u), \qquad f_2(u) = \Lambda_0(u),$$
$$h_1(z) = \exp(\beta^T z), \qquad h_2(z) = \exp(\beta^T z) - \exp(\beta_0^T z).$$

Then (5.1) implies that $\Lambda_0(u) \exp(\beta_0^T z) = \Lambda(u) \exp(\beta^T z)$ a.e. $\nu_1$, or, equivalently

$$0 = \{\Lambda(u) - \Lambda_0(u)\} e^{\beta^T z} + \Lambda_0(u)(e^{\beta^T z} - e^{\beta_0^T z})$$
$$= f_1(u) h_1(z) + f_2(u) h_2(z) \qquad \text{a.e. } \nu_1.$$

Since $\mu_1 \times H$ is absolutely continuous with respect to $\nu_1$ by assumption C2, equality holds in the last display a.e. with respect to $\mu_1 \times H$. By multiplying across the identity in the last display by $ab$, integrating with respect to the measure $\mu_1 \times H$, and then applying Fubini's theorem, it follows that

$$0 = \int f_1 a \, d\mu_1 \int h_1 b \, dH + \int f_2 a \, d\mu_1 \int h_2 b \, dH$$

for all measurable functions $a = a(u)$ and $b = b(z)$. The choice of $a = f_1 1_A$ for $A \in \mathcal{B}_1$ and $b = h_1 1_B$ for $B \in \mathcal{B}_d$ yields

$$0 = \int f_1^2 1_A \, d\mu_1 \int h_1^2 1_B \, dH + \int f_1 f_2 1_A \, d\mu_1 \int h_1 h_2 1_B \, dH;$$

the choice of $a = f_2 1_A$ (for the same $A \in \mathcal{B}_1$) and $b = h_2 1_B$ (for the same set $B \in \mathcal{B}_d$) yields

$$0 = \int f_1 f_2 1_A \, d\mu_1 \int h_1 h_2 1_B \, dH + \int f_2^2 1_A \, d\mu_1 \int h_2^2 1_B \, dH.$$

Thus we have

$$\int f_1^2 1_A \, d\mu_1 \int h_1^2 1_B \, dH = -\int f_1 f_2 1_A \, d\mu_1 \int h_1 h_2 1_B \, dH$$
$$= \int f_2^2 1_A \, d\mu_1 \int h_2^2 1_B \, dH$$

for all $A \in \mathcal{B}_1$ and $B \in \mathcal{B}_d$. By Fubini's theorem, this yields

$$\int_{A \times B} f_1^2 h_1^2 \, d(\mu_1 \times H) = \int_{A \times B} f_2^2 h_2^2 \, d(\mu_1 \times H)$$

for all such sets $A, B$. But this implies that the measures $\gamma_1$ and $\gamma_2$ defined by $\gamma_j(A \times B) = \int_{A \times B} f_j^2 h_j^2 \, d(\mu_1 \times H)$, $j = 1, 2$, are equal for all the product sets $A \times B$, and hence, by a standard monotone class argument, we conclude that $\gamma_1 = \gamma_2$ as measures on $([0, \tau] \times \mathbb{R}^d, \mathcal{B}_1[0, \tau] \times \mathcal{B}_d)$. It follows that $f_1^2(u) h_1^2(z) = f_2^2(u) h_2^2(z)$ a.e. with respect to $\mu_1 \times H$. Thus we conclude that

$$\frac{f_1^2(u)}{f_2^2(u)} = \frac{h_2^2(z)}{h_1^2(z)} \qquad \text{a.e. on } \{(u, z) : f_1^2(u) > 0, h_1^2(z) > 0\},$$

or, in other words,

$$\left(\frac{\Lambda(u)}{\Lambda_0(u)} - 1\right)^2 = (1 - \exp((\beta_0 - \beta)^T z))^2$$

a.e. with respect to $\mu_1 \times H$. This implies that (5.2) holds in view of C7. Integrating across this identity with respect to $\mu_1$ yields

$$\int \left(\frac{\Lambda(u)}{\Lambda_0(u)} - 1\right)^2 d\mu_1(u) = (1 - \exp((\beta_0 - \beta)^T z))^2 \mu_1([0, \tau]) \qquad \text{a.e. } H,$$

and hence the right-hand side is a constant a.e. $H$. But this implies that $\beta = \beta_0$ in view of C7. Combining this with the last display shows that (5.2) holds.

For any given $\varepsilon > 0$, let $\tilde{\theta}_n^{ps} = (\hat{\beta}_n^{ps}, (1-\varepsilon)\hat{\Lambda}_n^{ps} + \varepsilon \Lambda_0) = (\hat{\beta}_n^{ps}, \hat{\Lambda}_n^{ps}) + \varepsilon(0, \Lambda_0 - \hat{\Lambda}_n^{ps})$. Since $\mathbb{M}_n^{ps}(\hat{\theta}_n^{ps}) \geq \mathbb{M}_n^{ps}(\tilde{\theta}_n^{ps}) = \mathbb{M}_n^{ps}(\hat{\theta}_n^{ps} + \varepsilon(0, \Lambda_0 - \hat{\Lambda}_n^{ps}))$, it follows that

$$0 \geq \lim_{\varepsilon \downarrow 0} \frac{\mathbb{M}_n^{ps}(\hat{\theta}_n^{ps} + \varepsilon(0, \Lambda_0 - \hat{\Lambda}_n^{ps})) - \mathbb{M}_n^{ps}(\hat{\theta}_n^{ps})}{\varepsilon}$$
$$= \mathbb{P}_n \left[\sum_{j=1}^{K} \left\{\frac{\mathbb{N}_{Kj}}{\hat{\Lambda}_{nKj}^{ps} \exp(\hat{\beta}_n^{psT} Z)} - 1\right\} (\Lambda_{0Kj} - \hat{\Lambda}_{nKj}^{ps}) \exp(\hat{\beta}_n^{psT} Z)\right],$$

where $\hat{\Lambda}^{ps}_{nKj} = \hat{\Lambda}^{ps}_n(T_{K,j})$. This yields

$$\mathbb{P}_n\left[\sum_{j=1}^K\left\{\mathbb{N}_{Kj}\frac{\Lambda_{0Kj}}{\hat{\Lambda}^{ps}_{nKj}} + \hat{\Lambda}^{ps}_{nKj}\exp(\hat{\beta}^{psT}_n Z)\right\}\right]$$

$$\leq \mathbb{P}_n\left\{\sum_{j=1}^K (\mathbb{N}_{Kj} + \Lambda_{0Kj}\exp(\hat{\beta}^{psT}_n Z))\right\}$$

$$\leq C\mathbb{P}_n\left\{\sum_{j=1}^K (\mathbb{N}_{K,j} + \Lambda_0(T_{K,j}))\right\} \underset{\text{a.s.}}{\to} CP\left\{\sum_{j=1}^K (\mathbb{N}_{K,j} + \Lambda_0(T_{K,j}))\right\},$$

by C1–C3 and the strong law of large numbers. (Here $C$ represents a constant. In the sequel $C$ appearing in different lines may represent different constants.) The limit on the right-hand side is finite. On the other hand,

$$\limsup_{n\to\infty}\mathbb{P}_n\left[\sum_{j=1}^K\left\{\mathbb{N}_{Kj}\frac{\Lambda_{0Kj}}{\hat{\Lambda}^{ps}_{nKj}} + \hat{\Lambda}^{ps}_{nKj}\exp(\hat{\beta}^{psT}_n Z)\right\}\right]$$

$$\geq \limsup_{n\to\infty}\mathbb{P}_n\left\{\sum_{j=1}^K 1_{[b,\tau]}(T_{K,j})\hat{\Lambda}^{ps}_{nKj}\exp(\hat{\beta}^{psT}_n Z)\right\}$$

$$\geq C\limsup_{n\to\infty}\hat{\Lambda}^{ps}_n(b)\mathbb{P}_n\left(\sum_{j=1}^K 1_{[b,\tau]}(T_{K,j})\right)$$

$$= C\limsup_{n\to\infty}\hat{\Lambda}^{ps}_n(b)\mu_1([b,\tau]).$$

Hence $\hat{\Lambda}^{ps}_n(t)$ is uniformly bounded almost surely for $t \in [0,b]$ if $\mu_1([b,\tau]) > 0$ for some $0 < b < \tau$ or for $t \in [0,\tau]$ if $\mu_1(\{\tau\}) > 0$. By the Helly selection theorem and the compactness of $\mathcal{R} \times \mathcal{F}$, it follows that $\hat{\theta}^{ps}_n = (\hat{\beta}^{ps}_n, \hat{\Lambda}^{ps}_n)$ has a subsequence $\hat{\theta}^{ps}_{n'} = (\hat{\beta}^{ps}_{n'}, \hat{\Lambda}^{ps}_{n'})$ converging to $\theta^+ = (\beta^+, \Lambda^+)$, where $\Lambda^+$ is an increasing bounded function defined on $[0,b]$ for a $b < \tau$ and it can be defined on $[0,\tau]$ if $\mu_1(\{\tau\}) > 0$. Following the same argument as in proving Theorem 4.1 of Wellner and Zhang [21], we can show that $\mathbb{M}^{ps}(\theta^+) \geq \mathbb{M}^{ps}(\theta_0)$. Since $\mathbb{M}^{ps}(\theta_0) \geq \mathbb{M}^{ps}(\theta^+)$, by the argument above (5.1), we conclude that $\mathbb{M}^{ps}(\theta^+) = \mathbb{M}^{ps}(\theta_0)$. Then (5.2) implies that $\beta^+ = \beta_0$ and $\Lambda^+ = \Lambda_0$ a.e. in $\mu_1$. Finally, the dominated convergence theorem yields the strong consistency of $(\hat{\beta}^{ps}_n, \hat{\Lambda}^{ps}_n)$ in the metric $d_1$.

Now we turn to the maximum likelihood estimator. Let $\mathbb{M}_n(\theta) = n^{-1}l_n(\beta,\Lambda) = \mathbb{P}_n m_\theta(X)$ and $\mathbb{M}(\theta) = Pm_\theta(X)$, where

$$m_\theta(X) = \sum_{j=1}^K \{\Delta\mathbb{N}_{Kj}\log\Delta\Lambda_{Kj} + \Delta\mathbb{N}_{Kj}\beta^T Z - \Delta\Lambda_{Kj}\exp(\beta^T Z)\}.$$

Much as in the pseudo-likelihood case, $\mathbb{M}$ has $\theta_0 = (\beta_0, \Lambda_0)$ as its unique maximum point, and

$$\beta = \beta_0 \quad \text{and} \quad \Lambda(v) - \Lambda(u) = \Lambda_0(v) - \Lambda_0(u)$$

(5.3)

a.e. with respect to $\mu_2$.

The proof of consistency then proceeds along the same lines as for the pseudo-likelihood estimator; see Wellner and Zhang [22] for the detailed argument. The upshot is that $(\hat{\beta}_n, \hat{\Lambda}_n)$ is almost surely consistent in the metric $d_2$. $\square$

PROOF OF THEOREM 3.2. We derive the rate of convergence by checking the conditions in Theorem 3.2.5 of van der Vaart and Wellner [18]. Here we give a detailed proof for the first part of the theorem, and for the second we point out the differences in the proof from the first. Let

$$m_\theta^{ps}(X) = \sum_{j=1}^{K} \{\mathbb{N}_{Kj} \log \Lambda(T_{K,j}) + \mathbb{N}_{Kj} \beta^T Z - \Lambda(T_{K,j}) \exp(\beta^T Z)\}$$

with $\mathbb{N}_{Kj} = \mathbb{N}(T_{K,j})$ and $\mathbb{M}^{ps}(\theta) = P m_\theta^{ps}(X)$. We have

$$\mathbb{M}^{ps}(\theta_0) - \mathbb{M}^{ps}(\theta)$$

$$= E_{(Z,K,\underline{T}_K)} \left[ \sum_{j=1}^{K} \Lambda(T_{K,j}) \exp(\beta^T Z) h\left\{ \frac{\Lambda_0(T_{K,j}) \exp(\beta_0^T Z)}{\Lambda(T_{K,j}) \exp(\beta^T Z)} \right\} \right].$$

Since $h(x) \geq (1/4)(x-1)^2$ for $0 \leq x \leq 5$, for any $\theta$ in a sufficiently small neighborhood of $\theta_0$

$$\mathbb{M}^{ps}(\theta_0) - \mathbb{M}^{ps}(\theta)$$

(5.4) $$\geq \frac{1}{4} E_{(Z,K,\underline{T}_K)} \left[ \sum_{j=1}^{K} \Lambda(T_{K,j}) \exp(\beta^T Z) \left\{ \frac{\Lambda_0(T_{K,j}) \exp(\beta_0^T Z)}{\Lambda(T_{K,j}) \exp(\beta^T Z)} - 1 \right\}^2 \right]$$

$$\geq C \int \{\Lambda(u) e^{\beta^T z} - \Lambda_0(u) e^{\beta_0^T z}\}^2 \, d\nu_1(u,z)$$

by C1, C2 and C6.

Let $g(t) = \Lambda_t(U) \exp(\beta_t^T Z)$ with $\Lambda_t = t\Lambda + (1-t)\Lambda_0$ and $\beta_t = t\beta + (1-t)\beta_0$ for $0 \leq t \leq 1$ with $(U, Z) \sim \nu_1$. Then $\Lambda(U) \exp(\beta^T Z) - \Lambda_0(U) \exp(\beta_0^T Z) = g(1) - g(0)$ and hence, by the mean value theorem, there exists a $0 \leq \xi \leq 1$ such that $g(1) - g(0) = g'(\xi)$. Since

$$g'(\xi) = \exp(\beta_\xi^T Z)[(\Lambda - \Lambda_0)(U) + \{\Lambda_0 + \xi(\Lambda - \Lambda_0)\}(U)(\beta - \beta_0)^T Z]$$

$$= \exp(\beta_\xi^T Z)[(\Lambda - \Lambda_0)(U)\{1 + \xi(\beta - \beta_0)^T Z\} + (\beta - \beta_0)^T Z \Lambda_0(U)],$$

from (5.4) we have

$$P\{m^{ps}_{\theta_0}(X) - m^{ps}_\theta(X)\}$$
$$\geq C \int [(\Lambda - \Lambda_0)(u)\{1 + \xi(\beta - \beta_0)^T z\} + (\beta - \beta_0)^T z \Lambda_0(u)]^2 \, d\nu_1(u, z)$$
$$= \nu_1\{g_1 h + g_2\}^2,$$

where $g_1(U, Z) \equiv (\beta - \beta_0)^T Z \Lambda_0(U)$, $g_2(U) = (\Lambda - \Lambda_0)(U)$ and $h(U, Z) = 1 + \xi(\Lambda - \Lambda_0)(U)/\Lambda_0(U)$ in the notation of Lemma 8.8, page 432, van der Vaart [17]. To apply van der Vaart's lemma we need to bound $[\nu_1(g_1 g_2)]^2$ by a constant less than one times $\nu_1(g_1^2)\nu_1(g_2^2)$. For the moment we write expectations under $\nu_1$ as $E_1$. But by the Cauchy–Schwarz inequality and then computing conditionally on $U$ we have

$$[E_1(g_1 g_2)]^2$$
$$= \{E_1[E_1(g_1 g_2|U)]\}^2 \leq E_1\{g_2^2\}E_1\{[E_1(g_1|U)]^2\}$$
$$= E_1\{g_2^2\}E_1\{\Lambda_0^2(U)[E_1((\beta - \beta_0)^T Z|U)]^2\}$$
$$= E_1\{g_2^2\}E_1\{\Lambda_0^2(U)E_1[(\beta - \beta_0)^T (Z - (Z - E_1(Z|U)))^{\otimes 2}(\beta - \beta_0)|U]\}$$
$$\leq (1 - \eta)E_1\{g_2^2\}E_1\{\Lambda_0^2(U)(\beta - \beta_0)^T E_1(ZZ^T|U)(\beta - \beta_0)^T\}$$
$$= (1 - \eta)E_1\{g_2^2\}E_1\{g_1^2\},$$

where the last inequality follows from C13. By van der Vaart's lemma,

$$\nu_1\{g_1 h + g_2\}^2 \geq C\{\nu_1(g_1^2) + \nu_1(g_2^2)\}$$
$$= C\{|\beta - \beta_0|^2 + \|\Lambda - \Lambda_0\|^2_{L_2(\mu_1)}\} = Cd_1^2(\theta, \theta_0).$$

To derive the rate of convergence, next we need to find a $\phi_n(\sigma)$ such that

$$E \sup_{d_1(\theta, \theta_0) < \sigma} |\mathbb{G}_n(m^{ps}_\theta(X) - m^{ps}_{\theta_0}(X))| \leq C\phi_n(\sigma).$$

We let $\mathcal{M}^1_\delta(\theta_0) = \{m^{ps}_\theta(X) - m^{ps}_{\theta_0}(X) : d_1(\theta, \theta_0) < \delta\}$ be the class of differences. We shall find an upper bound for the bracketing entropy numbers of this class. We also let $\mathcal{F}_\delta = \{\Lambda \in \mathcal{F} : \|\Lambda - \Lambda_0\|_{L_2(\mu_1)} \leq \delta\}$. Since $\mathcal{F}_\delta$ is a class of monotone nondecreasing functions, by Theorem 2.7.5 of van der Vaart and Wellner [18], for any $\varepsilon > 0$, there exists a set of brackets $[\Lambda_1^l, \Lambda_1^r], [\Lambda_2^l, \Lambda_2^r], \ldots, [\Lambda_q^l, \Lambda_q^r]$ with $q \leq \exp(M/\varepsilon)$, such that for any $\Lambda \in \mathcal{F}_\delta$, $\Lambda_i^l(t) \leq \Lambda(t) \leq \Lambda_i^r(t)$ for all $t \in O[T]$ and some $1 \leq i \leq q$, and $\int \{\Lambda_i^r(u) - \Lambda_i^l(u)\}^2 \, d\mu_1(u) \leq \varepsilon^2$. (Here we use the fact that $\mu_1$ is a finite measure under our hypotheses, and hence can be normalized to be a probability measure.)

For sufficiently small $\varepsilon > 0$ and $\delta > 0$, we can construct the bracketing functions so that $\Lambda_i^r(t) - \Lambda_i^l(t) \leq \gamma_1$ and $\Lambda_i^l(t) \geq \gamma_2$ with $\gamma_1, \gamma_2 > 0$ for all $t \in$ $O[T]$ and $1 \leq i \leq q$. Here is the proof for this claim: For any $\Lambda \in \mathcal{F}_\delta$, the result of Lemma 7.1 implies that $\Lambda_0(t) - \varepsilon_1 \leq \Lambda(t) \leq \Lambda_0(t) + \varepsilon_1$ for a sufficiently small $\varepsilon_1 > 0$ [$\varepsilon_1$ can be chosen as $(\delta/C)^{2/3}$ in view of Lemma 7.1] and for all $t \in O[T]$. For any $1 \leq i \leq q$, there is a $\Lambda \in \mathcal{F}_\delta$ such that $\|\Lambda_i^r - \Lambda\|_{L_2(\mu_1)} \leq \varepsilon$ and $\|\Lambda - \Lambda_i^l\|_{L_2(\mu_1)} \leq \varepsilon$, which implies that $\|\Lambda_i^r - \Lambda_0\|_{L_2(\mu_1)} \leq \varepsilon^* (\varepsilon^* = \sqrt{\varepsilon^2 + \delta^2})$ and $\|\Lambda_i^l - \Lambda_0\|_{L_2(\mu_1)} \leq \varepsilon^*$. By Lemma 7.1, this yields that $\Lambda_i^r(t) \leq \Lambda_0(t) + \varepsilon_2$ and $\Lambda_i^l(t) \geq \Lambda_0(t) - \varepsilon_2$ for a sufficiently small $\varepsilon_2 > 0$. [$\varepsilon_2$ can be chosen as $(\varepsilon^*/C)^{2/3}$.] Therefore our claim is justified by letting $\gamma_1 = 2\varepsilon_2$ and $\gamma_2 = \Lambda_0(\sigma) - \varepsilon_2$, in view of C8.

Since $\beta \in \mathcal{R}$, a compact set in $\mathbb{R}^d$, we can construct an $\varepsilon$-net for $\mathcal{R}$, $\beta_1, \beta_2, \ldots, \beta_p$ with $p = [(M'/\varepsilon^d)]$, such that for any $\beta \in \mathcal{R}$, there is an $s$ such that

$$|\beta^T Z - \beta_s^T Z| \leq \varepsilon \quad \text{and} \quad |\exp(\beta^T Z) - \exp(\beta_s^T Z)| \leq C\varepsilon.$$

Therefore we can construct a set of brackets for $\mathcal{M}_\delta^1(\theta_0)$ as follows:

$$[m_{i,s}^{ps^l}(X), m_{i,s}^{ps^r}(X)], \quad \text{for } i = 1, 2, \ldots, q; s = 1, 2, \ldots, p,$$

where

$$m_{i,s}^{ps^l}(X) = \sum_{j=1}^{K} [\mathbb{N}_{Kj} \log \Lambda_i^l(T_{K,j}) + \mathbb{N}_{Kj}(\beta_s^T Z - \varepsilon)$$
$$- \Lambda_i^r(T_{K,j})\{\exp(\beta_s^T Z) + C\varepsilon\}] - m_{\theta_0}(X)$$

and

$$m_{i,s}^{ps^r}(X) = \sum_{j=1}^{K} [\mathbb{N}_{Kj} \log \Lambda_i^r(T_{K,j}) + \mathbb{N}_{Kj}(\beta_s^T Z + \varepsilon)$$
$$- \Lambda_i^l(T_{K,j})\{\exp(\beta_s^T Z) - C\varepsilon\}] - m_{\theta_0}(X).$$

In what follows, we show that $\|f_{i,s}(X)\|_{P,B}^2 = \|m_{i,s}^{ps^r}(X) - m_{i,s}^{ps^l}(X)\|_{P,B}^2 \leq C\varepsilon^2$, where $\|\cdot\|_{P,B}$ is the "Bernstein norm" defined by $\|f\|_{P,B} = \{2P(e^{|f|} - 1 - |f|)\}^{1/2}$ (see van der Vaart and Wellner [18], page 324). Since $2(e^x - 1 - x) \leq x^2 e^x$ for $x \geq 0$, it follows that $\|f\|_{P,B}^2 \leq P(e^{|f|}|f|^2)$. Therefore, $\|f_{i,s}(X)\|_{P,B}^2 \leq P(e^{|f_{i,s}(X)|}|f_{i,s}(X)|^2)$. By writing out $f_{i,s}(X) = m_{i,s}^{ps^r}(X) - m_{i,s}^{ps^l}(X)$, we find that

$$|f_{i,s}(X)| \leq \mathbb{N}_{KK} \sum_{j=1}^{K} |(\log \Lambda_i^r(T_{K,j}) - \log \Lambda_i^l(T_{K,j}) + 2\varepsilon)|$$
$$+ \exp(\beta_s^T Z) \sum_{j=1}^{K} (\Lambda_i^r(T_{K,j}) - \Lambda_i^l(T_{K,j}))$$

$$+ C\varepsilon \sum_{j=1}^{K}(\Lambda_i^r(T_{K,j}) + \Lambda_i^l(T_{K,j})).$$

Since

(5.5)
$$\log y = \log x + (x + \xi(y-x))^{-1}(y-x)$$
$$\text{for } 0 < x \leq y, \text{ some } \xi \in [0,1],$$

we find that

$$\log \Lambda_i^r(T_{K,j}) \leq \log \Lambda_i^l(T_{K,j}) + \gamma_2^{-1}(\Lambda_i^r(T_{K,j}) - \Lambda_i^l(T_{K,j}))$$

by construction of $\Lambda_i^l$. Hence, by C9 and our claim above, we conclude further that $\sum_{j=1}^{K}(|\log \Lambda_i^r(T_{K,j}) - \log \Lambda_i^l(T_{K,j})| + 2\varepsilon)$, $\sum_{j=1}^{K}(\Lambda_i^r(T_{K,j}) - \Lambda_i^l(T_{K,j}))$ and $\sum_{j=1}^{K}(\Lambda_i^r(T_{K,j}) + \Lambda_i^l(T_{K,j}))$ are all uniformly bounded in $O[T]$. More explicitly, taking $\varepsilon_2 \leq 2^{-1}\Lambda_0(\sigma)$, noting that this implies $\delta \leq C(2^{-1}\Lambda_0(\sigma))^{3/2} \equiv \delta_0^{ps}$ with $C = (c_0/(24f_0))^{1/2}$ by Lemma 7.1, and using the relations $\varepsilon_2 = (\varepsilon^*/C)^{2/3}$, $\varepsilon^* = (\varepsilon^2 + \delta^2)^{1/2} \geq \delta$ and $\varepsilon_2 \leq 2^{-1}\Lambda_0(\sigma)$, we find that

$$\sum_{j=1}^{K}(\log \Lambda_i^r(T_{K,j}) - \log \Lambda_i^l(T_{K,j}) + 2\varepsilon)^2$$
$$\leq \sum_{j=1}^{K}\left(\frac{2\varepsilon_2}{\Lambda_0(\sigma) - \varepsilon_2} + 2\varepsilon\right)^2 \leq 4k_0(1 + \delta_0^{ps})^2.$$

Therefore, by arguing conditionally on $(Z, K, T_K)$ and using C10,

$$\|f_{i,s}(X)\|_{P,B}^2 \leq P(e^{|f_{i,s}(X)|}|f_{i,s}(X)|^2)$$
$$\leq CP\left\{e^{v\mathbb{N}_{KK}}\left[\mathbb{N}_{KK}^2 \sum_{j=1}^{K}(\log \Lambda_i^r(T_{K,j}) - \log \Lambda_i^l(T_{K,j}) + 2\varepsilon)^2\right.\right.$$
$$\left.\left. + \exp(2\beta_s^T Z)\sum_{j=1}^{K}(\Lambda_i^r(T_{K,j}) - \Lambda_i^l(T_{K,j}))^2 + C\varepsilon^2\right]\right\}.$$

By C6, C10 and Taylor expansion for $\log \Lambda_i^r(T_{K,j})$ at $\Lambda_i^l(T_{K,j})$ as shown above, we have

$$\|f_{i,s}(X)\|_{P,B}^2 \leq C\left\{E_{(K,\underline{T}_K)}\left[\sum_{j=1}^{K}(\Lambda_i^r(T_{K,j}) - \Lambda_i^l(T_{K,j}))^2\right] + \varepsilon^2\right\} \leq C\varepsilon^2.$$

This shows that the total number of $\varepsilon$-brackets for $\mathcal{M}_\delta^1(\theta_0)$ will be of the order $(M/\varepsilon)^d e^{C(M'/\varepsilon)}$ and hence $\log N_{[]}(\varepsilon, \mathcal{M}_\delta^1(\theta_0), \|\cdot\|_{P,B}) \leq C(1/\varepsilon)$.

We can similarly verify that $P(f_\theta^{ps}(X))^2 \leq C\delta^2$ for any $f_\theta^{ps}(X) = m_\theta^{ps}(X) - m_{\theta_0}^{ps}(X) \in \mathcal{M}_\delta^1(\theta_0)$. Hence by Lemma 3.4.3 of van der Vaart and Wellner [18],

$$E_P^* \|\mathbb{G}_n\|_{\mathcal{M}_\delta^1(\theta_0)} \leq C\tilde{J}_{[]}(\delta, \mathcal{M}_\delta^1(\theta_0), \|\cdot\|_{P,B})\left[1 + \frac{\tilde{J}_{[]}(\delta, \mathcal{M}_\delta^1(\theta_0), \|\cdot\|_{P,B})}{\delta^2\sqrt{n}}\right],$$

where

$$\tilde{J}_{[]}(\delta, \mathcal{M}_\delta^1(\theta_0), \|\cdot\|_{P,B}) = \int_0^\delta \sqrt{1 + \log N_{[]}(\varepsilon, \mathcal{M}_\delta^1(\theta_0), \|\cdot\|_{P,B})}\, d\varepsilon$$

$$= C\int_0^\delta \sqrt{1 + \frac{1}{\varepsilon}}\, d\varepsilon \leq C\int_0^\delta \varepsilon^{-1/2}\, d\varepsilon \leq C\delta^{1/2}.$$

Hence $\phi_n(\delta) = \delta^{1/2}(1 + \delta^{1/2}/(\delta^2\sqrt{n})) = \delta^{1/2} + \delta^{-1}/\sqrt{n}$. Then it is easy to see that $\phi_n(\delta)/\delta$ is a decreasing function of $\delta$, and $n^{2/3}\phi_n(n^{-1/3}) = n^{2/3}(n^{-1/6} + n^{1/3}n^{-1/2}) = 2\sqrt{n}$. So it follows by Theorem 3.2.5 of van der Vaart and Wellner [18] that

$$n^{1/3}d_1((\hat{\beta}_n^{ps}, \hat{\Lambda}_n^{ps}), (\beta_0, \Lambda_0)) = O_P(1).$$

For the maximum likelihood estimator $(\hat{\beta}_n, \hat{\Lambda}_n)$ the proof of the rate of convergence result as stated in Theorem 3.2 proceeds along the same lines as the rate result for the maximum pseudo-likelihood estimator given above, but with $g_1(U, V, Z) \equiv (\beta - \beta_0)^T Z \Delta\Lambda_0(U, V)$, $g_2(U, V) = (\Delta\Lambda - \Delta\Lambda_0)(U, V)$ and $h(U, V, Z) = 1 + \xi(\Delta\Lambda - \Delta\Lambda_0)(U, V)/\Delta\Lambda_0(U, V)$ in the application of van der Vaart's Lemma 8.8. For details see Wellner and Zhang [22]. $\square$

PROOF OF THEOREM 3.3. We give a detailed proof for the first part of the theorem, and only outline the differences in the proof for the second. We prove the theorem by checking the conditions A1–A6 of Theorem 6.1. Note that A1 holds with $\gamma = 1/3$ because of the rate of convergence given in Theorem 3.2. The criterion function with only one observation is given by $m^{ps}(\beta, \Lambda; X) = \sum_{j=1}^K \{\mathbb{N}_{Kj} \log \Lambda_{Kj} + \mathbb{N}_{Kj}\beta^T Z - e^{\beta^T Z}\Lambda_{Kj}\}$, and thus we have

$$m_1^{ps}(\beta, \Lambda; X) = \sum_{j=1}^K Z(\mathbb{N}_{Kj} - \Lambda(T_{K,j})\exp(\beta^T Z)),$$

$$m_2^{ps}(\beta, \Lambda; X)[h] = \sum_{j=1}^K \left(\frac{\mathbb{N}_{Kj}}{\Lambda_{Kj}} - \exp(\beta^T Z)\right)h_{Kj},$$

$$m_{11}^{ps}(\beta, \Lambda; X)[h] = -\sum_{j=1}^K \Lambda_{Kj} ZZ^T \exp(\beta^T Z),$$

$$m_{12}^{ps}(\beta, \Lambda; X)[h] = m_{21}^{ps'}(\beta, \Lambda; X)[h] = -\sum_{j=1}^K Z\exp(\beta^T Z)h_{Kj}$$

and

$$m_{22}^{ps}(\beta, \Lambda; X)[\mathbf{h}, h] = -\sum_{j=1}^{K} \frac{\mathbb{N}_{Kj}}{\Lambda_{Kj}^2} \mathbf{h}_{Kj} h_{Kj},$$

where $\Lambda_{Kj} = \Lambda(T_{K,j})$ and $h_{Kj} = \int_0^{T_{K,j}} h(t)\,d\Lambda(t)$ for $h \in L_2(\Lambda)$. A2 automatically holds by the model assumption (1.1). For A3, we need to find an $\mathbf{h}^*$ such that

$$\dot{S}_{12}^{ps}(\beta_0, \Lambda_0)[h] - \dot{S}_{22}^{ps}(\beta_0, \Lambda_0)[\mathbf{h}^*, h]$$
$$= P\{m_{12}^{ps}(\beta_0, \Lambda_0; X)[h] - m_{22}^{ps}(\beta_0, \Lambda_0; X)[\mathbf{h}^*, h]\} = 0,$$

for all $h \in L_2(\Lambda_0)$. Note that

$$P\{m_{12}^{ps}(\beta_0, \Lambda_0; X)[h] - m_{22}^{ps}(\beta_0, \Lambda_0; X)[\mathbf{h}^*, h]\}$$
$$= -E\left\{\sum_{j=1}^{K}\left[Ze^{\beta_0^T Z} - \frac{\mathbb{N}_{Kj}}{(\Lambda_{0Kj})^2}\mathbf{h}_{Kj}^*\right]h_{Kj}\right\}$$
$$= -E_{(K,\underline{T}_K,Z)}\left\{\sum_{j=1}^{K}\left[Ze^{\beta_0^T Z} - \frac{e^{\beta_0^T Z}\mathbf{h}_{Kj}^*}{\Lambda_{0Kj}}\right]h_{Kj}\right\}.$$

Therefore, an obvious choice of $\mathbf{h}^*$ is

$$\mathbf{h}_{Kj}^* = \Lambda_{0Kj} E(Ze^{\beta_0^T Z}|K, T_{K,j})/E(e^{\beta_0^T Z}|K, T_{K,j}) \equiv \Lambda_{0Kj} R^{ps}(K, T_{K,j}).$$

Hence

$$m^{*ps}(\beta_0, \Lambda_0; X) = m_1^{ps}(\beta_0, \Lambda_0; X) - m_2^{ps}(\beta_0, \Lambda_0; X)[\mathbf{h}^*]$$
$$= \sum_{j=1}^{K}\left\{Z(\mathbb{N}_{Kj} - e^{\beta_0^T Z}\Lambda_{0Kj})\right.$$
$$\left. - \left(\frac{\mathbb{N}_{Kj}}{\Lambda_{0Kj}} - e^{\beta_0^T Z}\right)\Lambda_{0Kj} R^{ps}(K, T_{K,j})\right\}$$
$$= \sum_{j=1}^{K}(\mathbb{N}_{Kj} - e^{\beta_0^T Z}\Lambda_{0Kj})[Z - R^{ps}(K, T_{K,j})],$$
$$A^{ps} = -\dot{S}_{11}^{ps}(\beta_0, \Lambda_0) + \dot{S}_{21}^{ps}(\beta_0, \Lambda_0)[\mathbf{h}^*]$$
$$= E\left\{\sum_{j=1}^{K}[\Lambda_{0Kj}e^{\beta_0^T Z}ZZ^T - e^{\beta_0^T Z}\Lambda_{0Kj}R^{ps}(K, T_{K,j})]\right\}$$
$$= E_{(K,\underline{T}_K,Z)}\left\{\sum_{j=1}^{K}\Lambda_{0Kj}e^{\beta_0^T Z}[Z - R^{ps}(K, T_{K,j})]Z^T\right\}$$

$$= E_{(K,\underline{T}_K,Z)}\left\{\sum_{j=1}^{K} \Lambda_{0Kj} e^{\beta_0^T Z}[Z - R^{ps}(K, T_{K,j})]^{\otimes 2}\right\}$$

and

$$B^{ps} = E m^{*ps}(\beta_0, \Lambda_0; X)^{\otimes 2}$$
$$= E_{(K,\underline{T}_K,Z)}\left\{\sum_{j,j'=1}^{K} C_{j,j'}^{ps}(Z)[Z - R^{ps}(K, T_{K,j})][Z - R^{ps}(K, T_{K,j})]^T\right\},$$

with

$$C_{j,j'}^{ps}(Z) = E[(\mathbb{N}_{Kj} - e^{\beta_0^T Z}\Lambda_{0Kj})(\mathbb{N}_{Kj'} - e^{\beta_0^T Z}\Lambda_{0Kj'})|Z, K, T_{K,j}, T_{K,j'}].$$

To verify A4, we note that the first part automatically holds, because

$$S_{1n}^{ps}(\hat{\beta}_n^{ps}, \hat{\Lambda}_n^{ps}) = \mathbb{P}_n m_1^{ps}(\hat{\beta}_n^{ps}, \hat{\Lambda}_n^{ps}; X) = 0$$

since $\hat{\beta}_n^{ps}$ satisfies the pseudo-score equation. Next we shall show that

$$S_{2n}^{ps}(\hat{\beta}_n^{ps}, \hat{\Lambda}_n^{ps})[\mathbf{h}^*]$$

(5.6)
$$= \mathbb{P}_n\left[\sum_{j=1}^{K} \frac{1}{\hat{\Lambda}_{nKj}^{ps}}\{\mathbb{N}_{Kj} - \hat{\Lambda}_{nKj}^{ps}\exp(\hat{\beta}_n^{psT}Z)\}\Lambda_{0Kj} R^{ps}(K, T_{K,j})\right]$$

$$= o_P(n^{-1/2})$$

with $\hat{\Lambda}_{nKj}^{ps} = \hat{\Lambda}_n^{ps}(T_{K,j})$. Since $(\hat{\beta}_n^{ps}, \hat{\Lambda}_n^{ps})$ maximizes $\mathbb{P}_n m_\theta^{ps}(X)$ over the feasible region, consider a path $\theta_\varepsilon = (\hat{\beta}_n^{ps}, \hat{\Lambda}_n^{ps} + \varepsilon h)$ for $h \in \mathcal{F}$. Then

$$\lim_{\varepsilon \downarrow 0} \frac{d}{d\varepsilon}\mathbb{P}_n m_{\theta_\varepsilon}^{ps}(X) = \mathbb{P}_n\left[\sum_{j=1}^{K} \frac{1}{\hat{\Lambda}_{nKj}^{ps}}\{\mathbb{N}_{Kj} - \hat{\Lambda}_{nKj}^{ps}\exp(\hat{\beta}_n^{psT}Z)\}h_{Kj}\right] = 0.$$

Now choose $h_{Kj} = \hat{\Lambda}_{nKj}^{ps} E(Z\exp(\beta_0^T Z)|K, T_{K,j})/E(\exp(\beta_0^T Z)|K, T_{K,j})$. Then to demonstrate (5.6), it suffices to show that

$$I = \mathbb{P}_n\left[\sum_{j=1}^{K} \frac{1}{\hat{\Lambda}_{nKj}^{ps}}\{\mathbb{N}_{Kj} - \hat{\Lambda}_{nKj}^{ps}\exp(\hat{\beta}_n^{psT}Z)\}(\Lambda_{0Kj} - \hat{\Lambda}_{nKj}^{ps})\alpha_{Kj}\right]$$

$$= o_P(n^{-1/2}),$$

where $\alpha_{Kj} = E(Z\exp(\beta_0^T Z)|K, T_{K,j})/E(\exp(\beta_0^T Z)|K, T_{K,j})$. But $I$ can be decomposed as $I = I_1 - I_2 + I_3$, where

$$I_1 = (\mathbb{P}_n - P)\left\{\sum_{j=1}^{K} \frac{\mathbb{N}_{Kj}}{\hat{\Lambda}_{nKj}^{ps}}(\Lambda_{0Kj} - \hat{\Lambda}_{nKj}^{ps})\alpha_{Kj}\right\},$$

$$I_2 = (\mathbb{P}_n - P)\left\{\sum_{j=1}^{K} \exp(\hat{\beta}_n^{psT}Z)(\Lambda_{0Kj} - \hat{\Lambda}_{nKj}^{ps})\alpha_{Kj}\right\}$$

and

$$I_3 = P\left[\sum_{j=1}^{K} \frac{1}{\hat{\Lambda}_{nKj}^{ps}}\{\mathbb{N}_{Kj} - \hat{\Lambda}_{nKj}^{ps}\exp(\hat{\beta}_n^{psT}Z)\}(\Lambda_{0Kj} - \hat{\Lambda}_{nKj}^{ps})\alpha_{Kj}\right].$$

We show that $I_1$, $I_2$ and $I_3$ are all $o_P(n^{-1/2})$. Let

$$\phi_1(X;\Lambda) = \sum_{j=1}^{K} \frac{\mathbb{N}_{Kj}}{\Lambda_{Kj}}(\Lambda_{0Kj} - \Lambda_{Kj})\alpha_{Kj},$$

$$\phi_2(X;\beta,\Lambda) = \sum_{j=1}^{K} \exp(\beta^T Z)(\Lambda_{0Kj} - \Lambda_{Kj})\alpha_{Kj},$$

and define two classes $\Phi_1(\eta)$ and $\Phi_2(\eta)$ as

$$\Phi_1(\eta) = \{\phi_1 : \Lambda \in \mathcal{F} \text{ and } \|\Lambda - \Lambda_0\|_{L_2(\mu_1)} \leq \eta\}$$

and

$$\Phi_2(\eta) = \{\phi_2 : (\beta,\Lambda) \in \mathcal{R} \times \mathcal{F} \text{ and } d_1((\beta,\Lambda),(\beta_0,\Lambda_0)) \leq \eta\}.$$

Using the same bracketing entropy arguments as used in deriving the rate of convergence, it follows that both $\Phi_1(\eta)$ and $\Phi_2(\eta)$ are $P$-Donsker classes under conditions C1, C6 and C8. Moreover, for the seminorm $\rho_P(f) = \{P(f - Pf)^2\}^{1/2}$, under conditions C1, C6, C8 and C9, we have $\sup_{\phi_1 \in \Phi_1(\eta)} \rho_P(\phi_1) \to 0$ and $\sup_{\phi_2 \in \Phi_2(\eta)} \rho_P(\phi_2) \to 0$ if $\eta \to 0$. Due to the relationship between $P$-Donsker and asymptotic equicontinuity (see Corollary 2.3.12 of van der Vaart and Wellner [18]), this yields $I_1 = o_P(n^{-1/2})$ and $I_2 = o_P(n^{-1/2})$. For $I_3$, we have

$$\begin{aligned}
I_3 &= P\left[\sum_{j=1}^{K} \frac{\{\mathbb{N}_{Kj} - \hat{\Lambda}_{nKj}^{ps}\exp(\hat{\beta}_n^{ps'}Z)\}}{\hat{\Lambda}_{nKj}^{ps}}(\Lambda_{0Kj} - \hat{\Lambda}_{nKj}^{ps})\alpha_{Kj}\right] \\
&= E\left[\sum_{j=1}^{K} \frac{\{\Lambda_{0Kj}\exp(\beta_0^T Z) - \hat{\Lambda}_{nKj}^{ps}\exp(\hat{\beta}_n^{psT}Z)\}}{\hat{\Lambda}_{nKj}^{ps}}(\Lambda_{0Kj} - \hat{\Lambda}_{nKj}^{ps})\alpha_{Kj}\right] \\
&= E\left[\sum_{j=1}^{K} \frac{\{(\Lambda_{0Kj} - \hat{\Lambda}_{nKj}^{ps})e^{\beta_0^T Z} + \hat{\Lambda}_{nKj}^{ps}(e^{\beta_0^T Z} - e^{\beta^T Z})\}}{\hat{\Lambda}_{nKj}^{ps}}(\Lambda_{0Kj} - \hat{\Lambda}_{nKj}^{ps})\alpha_{Kj}\right] \\
&\leq Cd_1^2\{(\hat{\beta}_n^{ps},\hat{\Lambda}_n^{ps}),(\beta_0,\Lambda_0)\},
\end{aligned}$$

by performing Taylor expansion of $\exp(\beta^T Z)$ at $\beta_0$ along with conditions C1, C6, C8 and the result of Lemma 7.1. Finally the rate of convergence yields $I_3 \leq Cn^{-2/3}$ in probability and thus $I_3 = o_P(n^{-1/2})$.

To verify A5, we note that

$$\sqrt{n}(S_{1n}^{ps} - S_1^{ps})(\beta, \Lambda) - \sqrt{n}(S_{1n}^{ps} - S_1^{ps})(\beta_0, \Lambda_0)$$
(5.7)
$$= \mathbb{G}_n \left[ \sum_{j=1}^{K} Z\{\Lambda_{0Kj} \exp(\beta_0^T Z) - \Lambda_{Kj} \exp(\beta^T Z)\} \right]$$

and

$$\sqrt{n}(S_{2n}^{ps} - S_2^{ps})(\beta, \Lambda)[\mathbf{h}^*] - \sqrt{n}(S_{2n}^{ps} - S_2^{ps})(\beta_0, \Lambda_0)[\mathbf{h}^*]$$
(5.8)
$$= \mathbb{G}_n \left( \sum_{j=1}^{K} \left[ \left( \frac{\mathbb{N}_{Kj}}{\Lambda_{Kj}} - \frac{\mathbb{N}_{Kj}}{\Lambda_{0Kj}} \right) - \{e^{\beta^T Z} - e^{\beta_0^T Z}\} \right] \Lambda_{0Kj} \alpha_{Kj} \right).$$

Let

$$a(\beta, \Lambda; X) = \sum_{j=1}^{K} Z\{\Lambda_{0Kj} \exp(\beta_0^T Z) - \Lambda_{Kj} \exp(\beta^T Z)\}$$

and

$$b(\beta, \Lambda; X) = \sum_{j=1}^{K} \left[ \mathbb{N}_{Kj} \left( \frac{1}{\Lambda_{Kj}} - \frac{1}{\Lambda_{0Kj}} \right) - \{\exp(\beta^T Z) - \exp(\beta_0^T Z)\} \right] \Lambda_{0Kj} \alpha_{Kj}.$$

For a $\eta > 0$, we define

$$A(\eta) = \{a(\beta, \Lambda; X) : d_1\{(\beta, \Lambda), (\beta_0, \Lambda_0)\} \leq \eta \text{ and } (\beta, \Lambda) \in \mathcal{R} \times \mathcal{F}\}$$

and

$$B(\eta) = \{b(\beta, \Lambda; X) : d_1\{(\beta, \Lambda), (\beta_0, \Lambda_0)\} \leq \eta \text{ and } (\beta, \Lambda) \in \mathcal{R} \times \mathcal{F}\}.$$

Then by applying the bracketing entropy arguments as in the rate of convergence proof, we can show that both $A(\eta)$ and $B(\eta)$ are $P$-Donsker classes under conditions C1, C6 and C8 and for a small enough $\eta > 0$. We can also show that $\sup_{a \in A(\eta)} \rho_P\{a(\beta, \Lambda; X)\} \to 0$ and $\sup_{b \in B(\eta)} \rho_P\{b(\beta, \Lambda; X)\} \to 0$ if $\eta \to 0$. Then the rate of convergence along with Corollary 2.3.12 of van der Vaart and Wellner [18] yields that

$$\sup_{|\beta - \beta_0| \leq \sigma_n, \|\Lambda - \Lambda_0\| \leq Cn^{-1/3}} |\mathbb{G}_n a(\beta, \Lambda; X)| = o_P(1)$$

and

$$\sup_{|\beta - \beta_0| \leq \sigma_n, \|\Lambda - \Lambda_0\| \leq Cn^{-1/3}} |\mathbb{G}_n b(\beta, \Lambda; X)| = o_P(1).$$

Hence A5 holds with $\gamma = 1/3$.

Finally, to verify A6, performing Taylor expansion of $m_1^{ps}(\beta, \Lambda; X)$ at the point $(\beta_0, \Lambda_0)$, we have

$$m_1^{ps}(\beta, \Lambda : X)$$
$$= \sum_{j=1}^{K} Z\{\mathbb{N}_{Kj} - \Lambda_{Kj} \exp(\beta^T Z)\}$$
$$= m_1^{ps}(\beta_0, \Lambda_0; X) + m_{11}^{ps}(\beta_0, \Lambda_0; X)(\beta - \beta_0) + m_{12}^{ps}(\beta_0, \Lambda_0; X)[\Lambda - \Lambda_0]$$
$$- \sum_{j=1}^{K} \exp(\beta_0^T Z) Z Z^T (\beta - \beta_0)(\Lambda_{Kj} - \Lambda_{0Kj})$$
$$- \tfrac{1}{2} \sum_{j=1}^{K} Z \exp(\beta_\xi^T Z) \Lambda_{0Kj} (\beta - \beta_0)^T Z Z^T (\beta - \beta_0),$$

where $\beta_\xi = \beta_0 + \xi(\beta - \beta_0)$ for some $0 < \xi < 1$. This yields

$$|S_1^{ps}(\beta, \Lambda) - S_1^{ps}(\beta_0, \Lambda_0) - \dot{S}_{11}^{ps}(\beta_0, \Lambda_0)(\beta - \beta_0) - \dot{S}_{12}^{ps}(\beta_0, \Lambda_0)[\Lambda - \Lambda_0]|$$
$$(5.9) \quad = \left| P\left\{ \sum_{j=1}^{K} \exp(\beta_0^T Z) Z Z^T (\beta - \beta_0)(\Lambda_{Kj} - \Lambda_{0Kj}) \right.\right.$$
$$\left.\left. + \tfrac{1}{2} \sum_{j=1}^{K} Z \exp(\beta_\xi^T Z) \Lambda_{0Kj} (\beta - \beta_0)^T Z Z^T (\beta - \beta_0) \right\} \right|.$$

Similarly, we have

$$|S_2^{ps}(\beta, \Lambda)[\mathbf{h}^*] - S_2^{ps}(\beta_0, \Lambda_0)[\mathbf{h}^*] - \dot{S}_{21}^{ps}(\beta_0, \Lambda_0)[\mathbf{h}^*](\beta - \beta_0)$$
$$- \dot{S}_{22}^{ps}(\beta_0, \Lambda_0)[\mathbf{h}^*, \Lambda - \Lambda_0]|$$
$$(5.10) \quad = \frac{1}{2}\left| P\left\{ \sum_{j=1}^{K} \frac{(\Lambda_{Kj} - \Lambda_{0Kj})^2}{\Lambda_{\zeta Kj}} \mathbb{N}_{Kj} \mathbf{h}^*_{kj} \right.\right.$$
$$\left.\left. + \sum_{j=1}^{K} Z \exp(\beta_\zeta^T Z)(\beta - \beta_0)^T Z Z^T (\beta - \beta_0) \right\} \right|,$$

where $\beta_\zeta = \beta_0 + \zeta(\beta - \beta_0)$ and $\Lambda_{\zeta Kj} = \Lambda_{0Kj} + \zeta(\Lambda_{Kj} - \Lambda_{0Kj})$ for some $0 \leq \zeta \leq 1$. Hence by C1, C3, C6, C7 and C8, it follows that (5.9) and (5.10) $\leq C\{|\beta - \beta_0|^2 + \|\Lambda - \Lambda_0\|_{L_2(\mu_1)}^2\}$, so A6 holds with $\alpha = 2$ and thus the proof for the first part of Theorem 3.3 is complete.

For the second part, first we note that with a single observation, $m(\beta, \Lambda; X) = \sum_{j=1}^{K}\{\Delta\mathbb{N}_{Kj} \log \Delta\Lambda_{Kj} + \Delta\mathbb{N}_{Kj}\beta^T Z - e^{\beta^T Z}\Delta\Lambda_{Kj}\}$, and hence

$$m_1(\beta, \Lambda; X) = \sum_{j=1}^{K} Z[\Delta\mathbb{N}_{Kj} - \Delta\Lambda_{Kj}e^{\beta^T Z}],$$

$$m_2(\beta, \Lambda; X)[h] = \sum_{j=1}^{K}\left[\frac{\Delta\mathbb{N}_{Kj}}{\Delta\Lambda_{Kj}} - e^{\beta^T Z}\right]\Delta h_{Kj},$$

$$m_{11}(\beta, \Lambda; X) = -\sum_{j=1}^{K}\Delta\Lambda_{Kj} Z Z^T e^{\beta^T Z},$$

$$m_{12}(\beta, \Lambda; X)[h] = m_{21}^T(\beta, \Lambda; X)[h] = -\sum_{j=1}^{K} Z e^{\beta^T Z}\Delta h_{Kj},$$

$$m_{22}(\beta, \Lambda; X)[\mathbf{h}, h] = -\sum_{j=1}^{K}\frac{\Delta\mathbb{N}_{Kj}}{(\Delta\Lambda_{Kj})^2}\Delta\mathbf{h}_{Kj}\Delta h_{Kj},$$

where $\Delta h_{Kj} = \int_{T_{K,j-1}}^{T_{K,j}} h\, d\Lambda$ for $h \in L_2(\Lambda)$. A1 holds with $\gamma = 1/3$ and the norm $\|\cdot\|$ being $L_2(\mu_2)$ because of the rate of convergence established in Theorem 4.2. A2 holds by the model specification (1.1). For A3, we need to find an $\mathbf{h}^*$ such that

$$\dot{S}_{12}(\beta_0, \Lambda_0)[h] - \dot{S}_{22}(\beta_0, \Lambda_0)[\mathbf{h}^*, h]$$
$$= P\{m_{12}(\beta_0, \Lambda_0; X)[h] - m_{22}(\beta_0, \Lambda_0; X)[\mathbf{h}^*, h]\} = 0,$$

for all $h \in L_2(\Lambda_0)$. Note that

$$P\{m_{12}(\beta_0, \Lambda_0; X)[h] - m_{22}(\beta_0, \Lambda_0; X)[\mathbf{h}^*, h]\}$$
$$= -E\left\{\sum_{j=1}^{K}\left[Z e^{\beta_0' Z} - \frac{\Delta\mathbb{N}_{Kj}}{(\Delta\Lambda_{0Kj})^2}\Delta\mathbf{h}_{Kj}^*\right]\Delta h_{Kj}\right\}$$
$$= -E_{(K,\underline{T}_K,Z)}\left\{\sum_{j=1}^{K}\left[Z e^{\beta_0' Z} - \frac{e^{\beta_0' Z}\Delta\mathbf{h}_{Kj}^*}{\Delta\Lambda_{0Kj}}\right]\Delta h_{Kj}\right\}.$$

Therefore, an obvious choice of $\mathbf{h}^*$ satisfies

$$\Delta\mathbf{h}_{Kj}^* = \Delta\Lambda_{0Kj} E(Z e^{\beta_0^T Z}|K, T_{K,j-1}, T_{K,j})/E(e^{\beta_0^T Z}|K, T_{K,j-1}, T_{K,j})$$
$$\equiv \Delta\Lambda_{0Kj} R(K, T_{K,j-1}, T_{K,j}).$$

Hence

$$m^*(\beta_0, \Lambda_0; X) = m_1(\beta_0, \Lambda_0; X) - m_2(\beta_0, \Lambda_0; X)[\mathbf{h}^*]$$

$$= \sum_{j=1}^{K} \bigg\{ Z(\Delta \mathbb{N}_{Kj} - e^{\beta_0^T Z} \Delta \Lambda_{0Kj})$$

$$- \bigg(\frac{\Delta \mathbb{N}_{Kj}}{\Delta \Lambda_{0Kj}} - e^{\beta_0^T Z}\bigg) \Delta \Lambda_{0Kj} R(K, T_{K,j-1} T_{K,j}) \bigg\}$$

$$= \sum_{j=1}^{K} (\Delta \mathbb{N}_{Kj} - e^{\beta_0^T Z} \Delta \Lambda_{0Kj})[Z - R(K, T_{K,j-1}, T_{K,j})],$$

$$A = -\dot{S}_{11}(\beta_0, \Lambda_0) + \dot{S}_{21}(\beta_0, \Lambda_0)[\mathbf{h}^*]$$

$$= E_{(K,\underline{T}_K,Z)} \bigg\{ \sum_{j=1}^{K} \Delta \Lambda_{0Kj} e^{\beta_0^T Z} [Z - R(K, T_{K,j-1}, T_{K,j})] Z^T \bigg\}$$

$$= E_{(K,\underline{T}_K,Z)} \bigg\{ \sum_{j=1}^{K} \Delta \Lambda_{0Kj} e^{\beta_0^T Z} [Z - R(K, T_{K,j-1}, T_{K,j})]^{\otimes 2} \bigg\}$$

and

$$B = E m^*(\beta_0, \Lambda_0; X)^{\otimes 2}$$

$$= E_{(K,\underline{T}_K,Z)} \bigg\{ \sum_{j,j'=1}^{K} C_{j,j'}(Z)[Z - R(K, T_{K,j-1}, T_{K,j})]^{\otimes 2} \bigg\}$$

with

$$C_{j,j'}(Z) = E[(\Delta \mathbb{N}_{Kj} - e^{\beta_0^T Z} \Delta \Lambda_{0Kj})(\Delta \mathbb{N}_{Kj'} - e^{\beta_0^T Z} \Delta \Lambda_{0Kj'})$$
$$|Z, K, T_{K,j-1}, T_{K,j}, T_{K,j'-1}, T_{K,j'}].$$

The rest of the proof of the maximum likelihood part of Theorem 3.3 parallels the proof for the pseudo-likelihood estimator; see Wellner and Zhang [22] for the details. □

**6. A general theorem on the asymptotic normality of semiparametric M-estimators.** In this section, we present a general theorem dealing with the asymptotic normality of semiparametric $M$-estimators of regression parameters when the rate of convergence of the estimator for nuisance parameters is slower than $n^{-1/2}$. We consider a general setting of a semiparametric model: given i.i.d. observations $X_1, X_2, \ldots, X_n$, we estimate unknown parameters $(\beta, \Lambda)$ by maximizing an objective function $n^{-1} \sum_{i=1}^{n} m(\beta, \Lambda; X_i) = \mathbb{P}_n m(\beta, \Lambda; X)$, where $\beta$ is a finite-dimensional parameter and $\Lambda$ is an infinite-dimensional parameter. If $m$ happens to be the log-likelihood function based on a single observation, then the estimator is simply the semiparametric maximum likelihood estimator. Our Theorem 6.1 here generalizes Theorem 6.1 of Huang [7] to accommodate the situation when the model is misspecified for the observed data. For a misspecified model, the information matrix calculated in Huang [7] is not relevant since (6.2) of Huang [7] is no longer valid. The notation we use here follows that of Huang [7].

Let $\beta = (\beta, \Lambda)$, where $\beta \in \mathbb{R}^d$ and $\Lambda$ is an infinite-dimensional parameter in the class $\mathcal{F}$. Suppose that $\Lambda_\eta$ is a parametric path in $\mathcal{F}$ through $\Lambda$, that is, $\Lambda_\eta \in \mathcal{F}$, and $\Lambda_\eta|_{\eta=0} = \Lambda$.

Let $\mathbf{H} = \{h : h = \frac{\partial \Lambda_\eta}{\partial \eta}|_{\eta=0}\}$ and for any $h \in \mathbf{H}$ we define

$$m_1(\beta, \Lambda; x) = \nabla_\beta m(\beta, \Lambda; x) \equiv \left(\frac{\partial m(\beta, \Lambda; x)}{\partial \beta_1}, \ldots, \frac{\partial m(\beta, \Lambda; x)}{\partial \beta_d}\right)^T,$$

$$m_2(\beta, \Lambda; x)[h] = \frac{\partial m(\beta, \Lambda_\eta; x)}{\partial \eta}\bigg|_{\eta=0},$$

$$m_{11}(\beta, \Lambda; x) = \nabla^2_\beta m(\beta, \Lambda; x),$$

$$m_{12}(\beta, \Lambda; x)[h] = \frac{\partial m_1(\beta, \Lambda_\eta; x)}{\partial \eta}\bigg|_{\eta=0},$$

$$m_{21}(\beta, \Lambda; x)[h] = \nabla_\beta m_2(\beta, \Lambda; x)[h],$$

$$m_{22}(\beta, \Lambda; x)[h_1, h_2] = \frac{\partial^2 m(\beta, \Lambda_{\eta_j}; x)}{\partial \eta_1 \partial \eta_2}\bigg|_{\eta_j=0, j=1,2}$$

$$\equiv \frac{\partial}{\partial \eta_2} m_2(\beta, \Lambda_{\eta_2}; x)[h_1]\bigg|_{\eta_2=0}.$$

We also define

$$S_1(\beta, \Lambda) = Pm_1(\beta, \Lambda; X), \qquad S_2(\beta, \Lambda)[h] = Pm_2(\beta, \Lambda; X)[h],$$
$$S_{1n}(\beta, \Lambda) = \mathbb{P}_n m_1(\beta, \Lambda; X), \qquad S_{2n}(\beta, \Lambda)[h] = \mathbb{P}_n m_2(\beta, \Lambda; X)[h],$$
$$\dot{S}_{11}(\beta, \Lambda) = Pm_{11}(\beta, \Lambda; X), \qquad \dot{S}_{22}(\beta, \Lambda)[h, h] = Pm_{22}(\beta, \Lambda; X)[h, h],$$

$$\dot{S}_{12}(\beta, \Lambda)[h] = \dot{S}_{21}^T(\beta, \Lambda)[h] = Pm_{12}(\beta, \Lambda; X)[h].$$

Furthermore, for $\mathbf{h} = (h_1, h_2, \ldots, h_d)^T \in \mathbf{H}^d$, where $h_j \in \mathbf{H}$ for $j = 1, 2, \ldots, d$, we denote

$$m_2(\beta, \Lambda; x)[\mathbf{h}] = (m_2(\beta, \Lambda; x)[h_1], \ldots, m_2(\beta, \Lambda; x)[h_d])^T,$$
$$m_{12}(\beta, \Lambda; x)[\mathbf{h}] = (m_{12}(\beta, \Lambda; x)[h_1], \ldots, m_{12}(\beta, \Lambda; x)[h_d]),$$
$$m_{21}(\beta, \Lambda; x)[\mathbf{h}] = (m_{21}(\beta, \Lambda; x)[h_1], \ldots, m_{21}(\beta, \Lambda; x)[h_d])^T,$$
$$m_{22}(\beta, \Lambda; x)[\mathbf{h}, h] = (m_{22}(\beta, \Lambda; x)[h_1, h], \ldots, m_{22}(\beta, \Lambda; x)[h_d, h])^T,$$

and define
$$S_2(\beta,\Lambda)[\mathbf{h}] = Pm_2(\beta,\Lambda;X)[\mathbf{h}], \qquad S_{2n}(\beta,\Lambda)[\mathbf{h}] = \mathbb{P}_n m_2(\beta,\Lambda;X)[\mathbf{h}],$$
$$\dot{S}_{12}(\beta,\Lambda)[\mathbf{h}] = Pm_{12}(\beta,\Lambda;X)[\mathbf{h}], \qquad \dot{S}_{21}(\beta,\Lambda)[\mathbf{h}] = Pm_{21}(\beta,\Lambda;X)[\mathbf{h}],$$
$$\dot{S}_{22}(\beta,\Lambda)[\mathbf{h},h] = Pm_{22}(\beta,\Lambda;X)[\mathbf{h},h].$$

To establish the asymptotic distribution for the $M$-estimator $\hat{\beta}_n$, we need the following assumptions:

A1. $|\hat{\beta}_n - \beta_0| = o_p(1)$ and $\|\hat{\Lambda}_n - \Lambda_0\| = O_p(n^{-\gamma})$ for some $\gamma > 0$ and some norm $\|\cdot\|$.

A2. $S_1(\beta_0,\Lambda_0) = 0$ and $S_2(\beta_0,\Lambda_0)[h] = 0$ for all $h \in \mathbf{H}$.

A3. There exists an $\mathbf{h}^* = (h_1^*,\ldots,h_d^*)^T$, where $h_j^* \in \mathbf{H}$ for $j = 1,\ldots,d$, such that
$$\dot{S}_{12}(\beta_0,\Lambda_0)[h] - \dot{S}_{22}(\beta_0,\Lambda_0)[\mathbf{h}^*,h] = 0,$$
for all $h \in \mathbf{H}$. Moreover, the matrix
$$A = -\dot{S}_{11}(\beta_0,\Lambda_0) + \dot{S}_{21}(\beta_0,\Lambda_0)[\mathbf{h}^*]$$
$$= -P(m_{11}(\beta_0,\Lambda_0;X) - m_{21}(\beta_0,\Lambda_0;X)[\mathbf{h}^*])$$
is nonsingular.

A4. The estimator $(\hat{\beta}_n, \hat{\Lambda}_n)$ satisfies
$$S_{1n}(\hat{\beta}_n,\hat{\Lambda}_n) = o_P(n^{-1/2}) \quad \text{and} \quad S_{2n}(\hat{\beta}_n,\hat{\Lambda}_n)[\mathbf{h}^*] = o_P(n^{-1/2}).$$

A5. For any $\delta_n \downarrow 0$ and $C > 0$
$$\sup_{|\beta-\beta_0|\leq\delta_n, \|\Lambda-\Lambda_0\|\leq Cn^{-\gamma}} |\sqrt{n}(S_{1n}-S_1)(\beta,\Lambda) - \sqrt{n}(S_{1n}-S_1)(\beta_0,\Lambda_0)|$$
$$= o_P(1)$$
and
$$\sup_{|\beta-\beta_0|\leq\delta_n, \|\Lambda-\Lambda_0\|\leq Cn^{-\gamma}} |\sqrt{n}(S_{2n}-S_2)(\beta,\Lambda)[\mathbf{h}^*]$$
$$- \sqrt{n}(S_{2n}-S_2)(\beta_0,\Lambda_0)[\mathbf{h}^*]| = o_P(1).$$

A6. For some $\alpha > 1$ satisfying $\alpha\gamma > 1/2$, and for $(\beta,\Lambda)$ in a neighborhood of $(\beta_0,\Lambda_0)$: $\{(\beta,\Lambda): |\beta-\beta_0|\leq\delta_n,\ \|\Lambda-\Lambda_0\|\leq Cn^{-\gamma}\}$,
$$|S_1(\beta,\Lambda) - S_1(\beta_0,\Lambda_0) - \dot{S}_{11}(\beta_0,\Lambda_0)(\beta-\beta_0) - \dot{S}_{12}(\beta_0,\Lambda_0)[\Lambda-\Lambda_0]|$$
$$= o(|\beta-\beta_0|) + O(\|\Lambda-\Lambda_0\|^\alpha)$$
and
$$|S_2(\beta,\Lambda)[\mathbf{h}^*] - S_2(\beta_0,\Lambda_0)[\mathbf{h}^*]$$
$$- \dot{S}_{21}(\beta_0,\Lambda_0)[\mathbf{h}^*](\beta-\beta_0) - \dot{S}_{22}(\beta_0,\Lambda_0)[\mathbf{h}^*,\Lambda-\Lambda_0]|$$
$$= o(|\beta-\beta_0|) + O(\|\Lambda-\Lambda_0\|^\alpha).$$

THEOREM 6.1. *Suppose that assumptions* A1–A6 *hold. Then*

$$\sqrt{n}(\hat{\beta}_n - \beta_0) = A^{-1}\sqrt{n}\mathbb{P}_n m^*(\beta_0, \Lambda_0; X) + o_{p^*}(1) \xrightarrow{d} N(0, A^{-1}B(A^{-1})^T),$$

*where* $m^*(\beta_0, \Lambda_0; x) = m_1(\beta_0, \Lambda_0; x) - m_2(\beta_0, \Lambda_0; x)[\mathbf{h}^*]$, $B = Em^*(\beta_0, \Lambda_0; X)^{\otimes 2} = E(m^*(\beta_0, \Lambda_0; X)m^*(\beta_0, \Lambda_0; X)^T)$, *and* $A$ *is given in assumption* A3.

PROOF. A1 and A5 yield

$$\sqrt{n}(S_{1n} - S_1)(\hat{\beta}_n, \hat{\Lambda}_n) - \sqrt{n}(S_{1n} - S_1)(\beta_0, \Lambda_0) = o_P(1).$$

Since $S_{1n}(\hat{\beta}_n, \hat{\Lambda}_n) = o_{p^*}(n^{-1/2})$ by A4 and $S_1(\beta_0, \Lambda_0) = 0$ by A2, it follows that

$$\sqrt{n}S_1(\hat{\beta}_n, \hat{\Lambda}_n) + \sqrt{n}S_{1n}(\beta_0, \Lambda_0) = o_P(1).$$

Similarly,

$$\sqrt{n}S_2(\hat{\beta}_n, \hat{\Lambda}_n) + \sqrt{n}S_{2n}(\beta_0, \Lambda_0)[\mathbf{h}^*] = o_P(1).$$

Combining these equalities and A6 yields

$$(6.1) \quad \begin{aligned} &\dot{S}_{11}(\beta_0, \Lambda_0)(\hat{\beta}_n - \beta_0) + \dot{S}_{12}(\beta_0, \Lambda_0)[\hat{\Lambda}_n - \Lambda_0] + S_{1n}(\beta_0, \Lambda_0) \\ &\quad + o(|\hat{\beta}_n - \beta_0|) + O(\|\hat{\Lambda}_n - \Lambda_0\|^\alpha) = o_P(n^{-1/2}) \end{aligned}$$

and

$$(6.2) \quad \begin{aligned} &\dot{S}_{21}(\beta_0, \Lambda_0)[\mathbf{h}^*](\hat{\beta}_n - \beta_0) + \dot{S}_{22}(\beta_0, \Lambda_0)[\mathbf{h}^*, \hat{\Lambda}_n - \Lambda_0] \\ &\quad + S_{2n}(\beta_0, \Lambda_0)[\mathbf{h}^*] + o(|\hat{\beta}_n - \beta_0|) + O(\|\hat{\Lambda}_n - \Lambda_0\|^\alpha) \\ &= o_P(n^{-1/2}). \end{aligned}$$

Because $\alpha\gamma > 1/2$, the rate of convergence assumption A1 implies $\sqrt{n}O(\|\hat{\Lambda}_n - \Lambda_0\|^\alpha) = o_P(1)$. Thus by A4 and (6.1) minus (6.2), it follows that

$$(\dot{S}_{11}(\beta_0, \Lambda_0) - \dot{S}_{21}(\beta_0, \Lambda_0)[\mathbf{h}^*])(\hat{\beta}_n - \beta_0) + o(|\hat{\beta}_n - \beta_0|)$$
$$= -(S_{1n}(\beta_0, \Lambda_0) - S_{2n}(\beta_0, \Lambda_0)[\mathbf{h}^*]) + o_P(n^{-1/2}),$$

that is,

$$-(A + o(1))(\hat{\beta}_n - \beta_0) = -\mathbb{P}_n m^*(\beta_0, \Lambda_0; X) + o_P(n^{-1/2}).$$

This yields

$$\begin{aligned}\sqrt{n}(\hat{\beta}_n - \beta_0) &= (A + o(1))^{-1}\sqrt{n}\mathbb{P}_n m^*(\beta_0, \Lambda_0; X) + o_P(1) \\ &\xrightarrow{d} N(0, A^{-1}B(A^{-1})^T).\end{aligned} \qquad \square$$

## 7. A technical lemma.

LEMMA 7.1. *Suppose that conditions* C8, C11 *and* C12 *hold, and that* $\Lambda \in \mathcal{F}$ *satisfies* $\|\Lambda - \Lambda_0\|_{L_2(\mu_1)} \leq \eta$. *Then there exists a constant* $C$ *independent of* $\Lambda$ *such that*

$$\sup_{t \in O[T]} |\Lambda(t) - \Lambda_0(t)| \leq (\eta/C)^{2/3}.$$

PROOF. Suppose that $t_0 \in O[T]$ satisfies

$$|\Lambda(t_0) - \Lambda_0(t_0)| \geq (1/2) \sup_{t \in O[T]} |\Lambda(t) - \Lambda_0(t)| \equiv \xi/2.$$

Then either $\Lambda(t_0) \geq \Lambda_0(t_0) + \xi/2$, or $\Lambda_0(t_0) \geq \Lambda(t_0) + \xi/2$; that is, $\Lambda(t_0) \leq \Lambda_0(t_0) - \xi/2$. In the first case we have

$$\begin{aligned}
\eta^2 &\geq \int \{\Lambda(t) - \Lambda_0(t)\}^2 \, d\mu_1(t) \\
&\geq \int_{t_0}^{\Lambda_0^{-1}(\xi/2 + \Lambda_0(t_0))} \{\Lambda(t) - \Lambda_0(t)\}^2 \dot\mu_1(t) \, dt \\
&\geq \int_{t_0}^{\Lambda_0^{-1}(\xi/2 + \Lambda_0(t_0))} \{\Lambda_0(t_0) + \xi/2 - \Lambda_0(t)\}^2 \dot\mu_1(t) \, dt \\
&\geq \int_{\Lambda_0(t_0)}^{\xi/2 + \Lambda_0(t_0)} \{\Lambda_0(t_0) + \xi/2 - x\}^2 \dot\mu_1(\Lambda_0^{-1}(x)) \frac{1}{\Lambda_0'\{\Lambda_0^{-1}(x)\}} \, dx \\
&\geq (c_0/f_0) \int_{\Lambda_0(t_0)}^{\xi/2 + \Lambda_0(t_0)} \{\Lambda_0(t_0) + \xi/2 - x\}^2 \, dx \geq \frac{c_0}{24 f_0} \xi^3 \\
&= \frac{c_0}{24 f_0} \left( \sup_{t \in O[T]} |\Lambda(t) - \Lambda_0(t)| \right)^3.
\end{aligned}$$

This yields the stated conclusion with $C \equiv \sqrt{c_0/(24 f_0)}$. In the second case the same conclusion holds by a similar argument. $\square$

The result of Lemma 7.1 can be extended to the interval $S[T] = (0, \tau)$ as long as C12 is valid on $S[T]$ and $\dot\mu_1(t)$ is uniformly bounded away from zero for $t \in S[T]$.

**Acknowledgments.** We owe thanks to an Associate Editor and two anonymous referees for their helpful and constructive comments and suggestions. We also thank Minggen Lu, a graduate student in Department of Biostatistics at the University of Iowa, for his computing support in the revision of this paper.

DEPARTMENT OF STATISTICS
UNIVERSITY OF WASHINGTON
BOX 354322
SEATTLE, WASHINGTON 98195-4322
USA
E-MAIL: jaw@stat.washington.edu

DEPARTMENT OF BIOSTATISTICS
UNIVERSITY OF IOWA
200 HAWKINS DRIVE, C22 GH
IOWA CITY, IOWA 52242-1009
USA
E-MAIL: ying-j-zhang@uiowa.edu